\documentclass[1pt,final,times]{preprint}
\usepackage{snapshot}
\usepackage[margin=3cm]{geometry}

\usepackage{framed,multirow}

\usepackage{latexsym}
\usepackage{amsmath,amsfonts,amssymb,amsthm,amsbsy,bm}		
\theoremstyle{theorem}
\theoremstyle{definition}\newtheorem{defn}{Definition}[section]
\theoremstyle{definition}\newtheorem{remark}{Remark}[section]
\usepackage[usenames,dvipsnames]{xcolor}
\usepackage{lipsum}

\usepackage{algorithmic}
\usepackage[T1]{fontenc}							
\usepackage[utf8]{inputenc}							

\usepackage[normalem]{ulem}
\usepackage{graphicx}								
\usepackage{microtype}								
\usepackage{hyperref}								
\usepackage{cleveref}								
\usepackage{booktabs}								
\usepackage{float}
\ifpdf
  \DeclareGraphicsExtensions{.eps,.pdf,.png,.jpg,.jpeg} 
\else
  \DeclareGraphicsExtensions{.eps}
\fi

\usepackage[colorinlistoftodos]{todonotes}

\newcommand{\diff}{\mathrm{d}}
\newcommand{\R}{\mathbb{R}}
\newcommand{\bu}{\boldsymbol{u}}

\newcommand{\bx}{\boldsymbol{x}}
\newcommand{\br}{\boldsymbol{r}}
\newcommand{\bnu}{\boldsymbol{\nu}}

\newcommand{\dr}{\,{\mathrm{d}\br}}
\newcommand{\ds}{\,{\mathrm{d}s}}
\newcommand{\btau}{\boldsymbol{\tau}}

\newcounter{todocounter}

\DeclareMathOperator*{\argmin}{{\rm arg\,min}}

\DeclareMathAlphabet{\pazocal}{OMS}{zplm}{m}{n}

\begin{document}
\begin{frontmatter}
  \title{Model hierarchies and higher-order discretisation of time-dependent thin-film free boundary problems with dynamic contact angle}
  \author{Dirk Peschka\fnref{myfootnote}}
  \address{Weierstrass Institute (WIAS), Mohrenstr. 39, 10117 Berlin, Germany}
  \author{Luca Heltai}
  \address{SISSA - International School for Advanced Studies, Via Bonomea 265, 34136 Trieste, Italy}
  \fntext[myfootnote]{corresponding author: peschka@wias-berlin.de}


  \begin{abstract}
    We present a mathematical and numerical framework for 
    thin-film fluid flows over planar surfaces including dynamic contact angles.
    In particular, we provide algorithmic details and an implementation of higher-order spatial and temporal discretisation of the underlying free boundary problem using the finite element method. The corresponding partial differential equation is based on a thermodynamically consistent energetic variational formulation of the problem using free energy and viscous dissipation in the bulk, on the surface, and at the moving contact line. 
    Model hierarchies for limits of strong and weak contact line dissipation are established, implemented and studied.
    We analyze the performance of the numerical algorithm and investigate the impact of the dynamic contact angle on the evolution of two benchmark problems: gravity-driven sliding droplets and the instability of a ridge. \\[-0.1em]
  \end{abstract}


\end{frontmatter}


\section{Introduction}
Flows of fluids with moving free surfaces and contact lines are ubiquitous in nature and applications, and manifest on different length scales \cite{de2013capillarity,bonn2009wetting,oron1997long}. 
While patterns of turbulent behavior tend to dominate fluid flows on large length-scales, on small length-scales surface tension plays a dominant role, and the resulting flows are mostly laminar.
It has been stated by Helmholtz, Rayleigh, and Korteweg \cite{helmholtz1868theorie,strutt1871some,korteweg1883xvii} a long time ago  and was shown by many since then 
\cite{Grmela1997,Oettinger1997,Morrison1998,Liu2001,Gay2009,Ishii2011,doi2011onsager,giga2018variational} that flows of viscous liquids have an underlying energetic variational structure. This structure applies to many if not all aspects of the thermodynamic and continuum description of the flow: the bulk flow and underlying rheological models, the thermodynamic description of the interface evolution as well as the description of moving contact lines\cite{shikhmurzaev1993moving,qian2006variational,ren2007boundary,ren2011derivation}. While for bulk flows the mathematical analysis, thermodynamical consistent modeling, and numerical treatment are well understood, the corresponding treatment of free boundaries -- and in particular of moving contact lines -- is still challenging and these aspects are often still disputed.

For applications it is often desirable to control a fluid in order to create or prevent the creation of certain flow patterns, \emph{e.g.}, \cite{huang2005spontaneous}, or of certain hydrodynamic instabilities \emph{e.g.}, \cite{sharma1996instability,seemann2001dewetting}. 
For instance, in dewetting flows one can steer the dynamics towards different stationary patterns depending on the magnitude of dissipation at the interface \cite{peschka2019signatures}. 
In this work we provide an energetic variational approach and develop a numerical framework for the simulation of fluid flow over solid substrates, where the dissipation in the bulk phase, on the interface, and near the contact line can be tuned separately, and simulated efficiently on a variety of length scales.

\begin{figure}[hb] 
  \centering
  \includegraphics[width=0.55\textwidth]{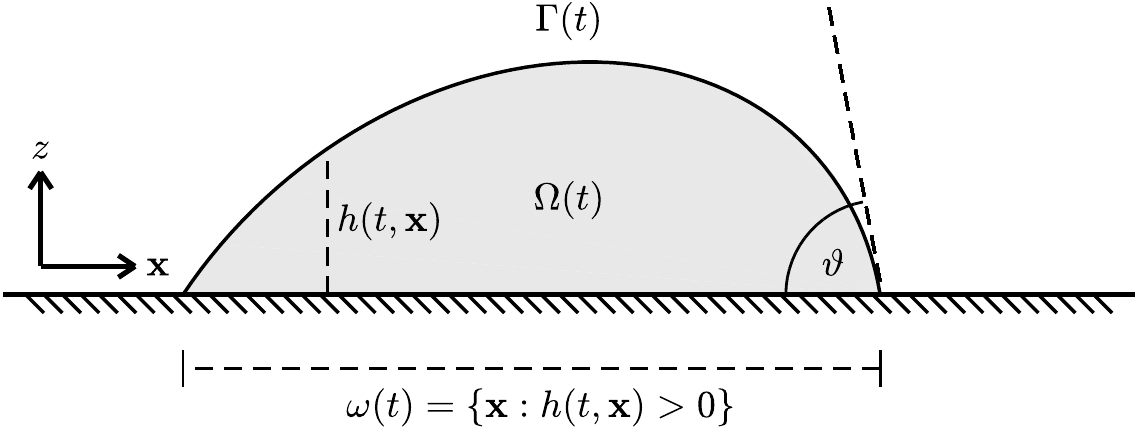}
  \caption{Sketch of geometry for free boundary problem based on $h(t):\omega(t)\to\mathbb{R}$.}
  \label{fig:cangle}
\end{figure}

%
Wetting and dewetting problems provide a physical setting which is relevant for applications and challenges our physical understanding of continuum fluid dynamics on various time- and length-scales. A key quantity of interest is the shape of the time-dependent domain occupied by a viscous fluid
\begin{subequations}
\label{eqn:geometricalshapes}
\begin{align}
  \Omega(t)=\Big\{\,\br=(\bx,z)\in\R^3:\bx=(x,y)\in\R^2,0<z<h(t,\bx)\,\Big\}\subset\R^3,
\end{align} 
which for convenience is often parametrized using a time-dependent \emph{height} function $h(t,\bx)\ge 0$. Assuming that the domain of interest can be parametrized in this way, we can additionally define
\begin{align}
\text{the wetted area} &\qquad\qquad \omega(t)=\{\bx\in\R^2:h(t,\bx)>0\},\\
\text{the liquid-vapor (gas) interface} &\qquad\qquad \Gamma(t)=\{(\bx,z)\in\R^3:\bx\in\omega(t),z=h(t,\bx)\},\\
\text{the contact line}&\qquad\qquad\partial\omega(t)\subset\R^2,
\end{align} 
\end{subequations}
as schematically presented in Fig.~\ref{fig:cangle}.

In the following we discuss the main implications of a continuum-mechanical description of the fluid flow in terms of a quadratic dissipation, following the spirit of the Rayleigh dissipation functional, \emph{e.g.} see also \cite{doi2011onsager}.
The material law for an incompressible Newtonian liquid with viscosity $\mu>0$ relates the flow field $\bu(t,\br)=\bigl(\bu_{\bx}(t,\br),u_z(t,\br)\bigr)\in\mathbb{R}^3$ at time $t$ in the fluid domain $\br\in\Omega(t)$ with the stress $\btau(\bu)(t,\br)\in\mathbb{R}^{3\times3}$ via
%
$\btau(\bu)=2\mu \mathbb{D}\bu$ 
 where
$\mathbb{D}\bu=\tfrac{1}{2}(\nabla\bu+\nabla\bu^\top)$.
%
While $\btau(\bu)$ is the prototypical stress tensor for bulk effects of viscous Newtonian flows, the no-slip condition
\begin{align}
\label{eqn:noslip}
\bu=\mathbf{0} \quad \text{on}\quad\omega(t)\times\{z=0\},
\end{align}
is the prototype boundary condition at solid-fluid interfaces $\omega(t)$.
In microfluidic settings it is the no-slip boundary condition \eqref{eqn:noslip} that needs to be extended in order
to allow for a
nonvanishing tangential velocity at liquid-solid interfaces. The established Navier-slip condition is
\begin{align}
\label{eqn:navierslip}
\bnu\cdot\bu=u_z=0,\qquad \mathbf{t}\cdot\btau(\bu) \bnu + 
\mu_\omega\mathbf{t}\cdot\bu=0\quad\text{on}\quad\omega(t)\times\{z=0\},
\end{align}
where $\bnu$ is the outer interface normal, $\mathbf{t}$ is a tangent vector, and the friction coefficient $\mu_\omega:=\mu/b>0$ is often associated with the so-called \emph{slip length} $b$. In particular, for flows near a moving contact line $\partial\omega$, the Navier-slip condition offers one out of various possible mechanism to remove the non-integrable stress singularity at a sliding contact line \cite{huh1971hydrodynamic}. Slip-lengths generally range from a few nanometers up to several micrometers, and are influenced by physical and chemical interface properties such as hydrophilicity, surface roughness, coating, or viscoelastic properties.
The validity of slip in microfluidic settings is undisputed nowadays \cite{lauga2007microfluidics}.
Different nonlinear generalisations of the Navier-slip condition have been proposed \cite{craig2001shear,thompson1997general}, which can be obtained by setting $\mu_\omega=\mu_\omega(|\bu|)$ in \eqref{eqn:navierslip}.
Thin-film models with different magnitudes of dissipation $\mu_\omega$ have been investigated in \cite{munch2005lubrication}. These models were extended to multiphase flows in \cite{jachalski2014impact}, where  the Navier-slip condition then reads
$\mathbf{t}\cdot\btau(\bu)\bnu+\mu_\omega\mathbf{t}\cdot[\![\bu]\!]_{\bnu}=0$, where $[\![\cdot]\!]_{\bnu}$ is the normal jump at the interface.


Free boundary problems become even more complex when considering moving contact lines. Contact lines are where solid, liquid, and vapor phases meet. 
Interfacial tensions between solid-vapor $\sigma_\text{SV}$, solid-liquid $\sigma_\text{SL}$ and liquid-vapor $\sigma$ lead to equilibrium contact angles $\vartheta=\vartheta_e$, which are determined by the famous Young-Dupr\'e equation
\begin{align}
\label{eqn:youngslaw}
\sigma\bigl(\cos(\vartheta_e)-1\bigr)=S \quad\text{ at }\quad \partial\omega(t)\times\{z=0\},
\end{align}
with spreading coefficient $S=\sigma_\text{SV}-(\sigma_\text{SL}+\sigma)$
on solid surfaces or by the Neumann triangle construction when considering soft elastic or liquid surfaces.
The situation becomes much more involved for  contact lines moving over physically realistic nonideal surfaces, \emph{e.g.}, which are not perfectly smooth, rigid, or chemically homogeneous.
Primarily for liquid-liquid-solid junctions, rolling motion was observed and thermodynamic descriptions have been developed, \emph{e.g.}, \cite{shikhmurzaev1997moving,bedeaux2004nonequilibrium}. In this work, however, we focus on sliding motion.
For sliding contact lines, models with contact angle hysteresis/stick-slip and empirical laws for contact line dynamics  are discussed in the literature, \emph{e.g.}, see \cite{snoeijer2013moving,bonn2009wetting}. In the Young-Dupr\'e equation \eqref{eqn:youngslaw} the spreading coefficient needs to be replaced using the Cassie-Baxter relation with the homogenized  coefficient $\bar{S}$ \cite{cassie1944wettability}.
However, the corresponding equilibrium angle might never be attained, since, due to pinning at microscopic defects, the contact line is stationary for contact angles $\vartheta_r < \vartheta < \vartheta_a$ and only starts to advance or recede beyond these values.

Models for contact line dynamics are often based on empirical laws and microscopic theories. For instance, assuming an energy-dissipation balance truncated at a microscopic length scale $\ell$ and with a macroscopic length scale $L$ away from the contact line, the theory by Cox and Voinov predicts
\begin{align*}
\vartheta^3-\vartheta_e^3
=\pm 9\,\text{Ca}\ln\bigl({L}/{\ell}\bigr),
\end{align*}
where the Capillary number $\text{Ca}=\mu U/\sigma$ depends on dynamic viscosity $\mu$, surface tension $\sigma$ and contact line velocity $U$ \cite{cox1986dynamics}.
Another contact line model by Ruschak and Hayes  produces a similar law, which, however, was derived from a molecular-kinetic model of wetting using dissipation at the contact line \cite{blake1973contact}.
Most models considered in the literature are of the form $W(\vartheta,\vartheta_e)=\pm\text{Ca}$, where usually $W$ is a polynomial with $W(\vartheta_e,\vartheta_e)=0$. From a thermodynamic point of view, the work by Ren \& E \cite{ren2007boundary} showed the motivation and mathematical modeling of nontrivial contact angle dynamics 
\begin{align}
\label{eqn:contactlinemodel}
	\sigma(\cos\vartheta_e-\cos\vartheta) = \mu_{\partial\omega} \bnu_x\cdot\bu_x,
\end{align}
for hydrodynamic continuum models, which itself was motivated by molecular kinetic models and points towards the fundamental character of this condition. 

Variational principles are used to model complex physical systems with partial differential equations derived from purely energetic arguments, see \emph{e.g.} \cite{otto2001geometry,peletier2014variational}. General variational approaches for complex fluids often rely on a classification of dissipative or nondissipative kinematics, \emph{e.g.}, see \cite{Grmela1997,beris1994thermodynamics,marsden2013introduction}.
In particular for flows of liquids over solids \cite{de1985wetting}, \emph{i.e.}, wetting and dewetting flows, such formulations amplify the physical and mathematical understanding of the corresponding free boundary problems.
The total dissipation by all these bulk, interface, and contact line contributions can be encoded it the functional $\Psi$, related to the \emph{Rayleighian} $\mathcal{R}(\bu)= \Psi(\bu)+\langle\mathrm{D}E_{\rm Stokes},\bu\rangle$, defined by
\begin{align}
  \label{eqn:newtondiss}
    2\Psi(\bu)=\langle K_\mathrm{Stokes}^{-1}(h)\bu,\bu\rangle=\int_\Omega 2\mu\mathbb{D}\bu:\nabla\bu\,\dr + \int_{\omega\times\{z=0\}} \mu_\omega |\bu|^2\ds + \int_{\partial\omega\times\{z=0\}} \mu_{\partial\omega} (\bu\cdot\bnu_x)^2\ds,
\end{align}
which is generated by a symmetric, positive operator $K_\mathrm{Stokes}(h)$ acting on the space of divergence-free velocities. This dissipation  encodes bulk-dissipation using the liquid viscosity $\mu>0$, a Navier-slip condition encoded in the friction coefficient $\mu_\omega\ge 0$, and a dynamic contact angle encoded in the friction coefficient $\mu_{\partial\omega}\ge 0$.


Assuming gravity and surface tension are the only contributing forces, the (potential/free) energy of the system is 
\begin{align}
\label{eqn:newtonenergy}
E_\mathrm{Stokes}(h)=\int_{\Omega(t)} \varrho (g_x\cdot\boldsymbol{x}+g_zz)\dr +  \int_{\Gamma(t)}\sigma\ds + \int_{\omega(t)\times\{z=0\}}(\sigma_\text{SL}-\sigma_\text{SV})\ds,
\end{align}
characterized by the geometrical shape of the domain $h$. The energy $E$ decreases during the evolution $\tfrac{\diff}{\diff t}E\bigl(h(t)\bigr)\le 0$. From now on we assume the evolution is a quasi-static in the sense that we can neglect the kinetic energy and there is an exact balance of change of potential energy and dissipation in the sense
\begin{align}
  \tfrac{\diff}{\diff t}E_\mathrm{Stokes}\bigl(h(t)\bigr) = \langle {\rm D}E_\mathrm{Stokes}(h),\bu\rangle=  -2\Psi(\bu)
  \equiv -\langle{\rm D}E_\mathrm{Stokes}(h),K_\mathrm{Stokes}(h){\rm D}E_\mathrm{Stokes}(h)\rangle\le 0,
\end{align}
where the operator $K_\mathrm{Stokes}:V^*\to V$ maps forces ${\rm D}E_\mathrm{Stokes}(h)\in V^*$ to velocities $\bu\in V$ and $\langle\cdot,\cdot\rangle:V^*\times V\to\R$ is the dual pairing. A formal but instructive calculation shows that the velocity can be found by solving
\begin{align}
\bu = -K_\mathrm{Stokes}(h){\rm D}E_\mathrm{Stokes}(h),
\end{align}
by solving the elliptic problem implictly defined through \eqref{eqn:newtondiss} with the right-hand-side being the shape-derivative
\begin{align}
\langle {\rm D}E_\mathrm{Stokes}(h),\mathbf{v}\rangle=\int_{\partial\Omega} \varrho \bigl(g_x\cdot\boldsymbol{x}+g_zz\bigr)\mathbf{v}\cdot\bnu\ds+\int_{\Gamma} \sigma \nabla_\Gamma \mathbf{v}\cdot\nabla_\Gamma\mathbf{r}\ds+\int_{\omega\times\{z=0\}}(\sigma_\text{SL}-\sigma_\text{SV})\nabla_\Gamma \mathbf{v}\cdot\nabla_\Gamma\mathbf{r}\ds,
\end{align}
i.e., the Fr\'echet derivative of the driving energy \eqref{eqn:newtonenergy} with respect to changes of the domain with the velocity $\mathbf{v}$.  Using incompressibility $\nabla_{\bx}\cdot\bu_{\bx}+\partial_z u_z=0$ and the kinematic condition $\partial_t h=u_z-\bu_{\bx}\cdot\nabla_{\bx}h$ one can show
\begin{align}
\partial_t h(t,\bx) = -\nabla_{\bx}\cdot\int_0^{h(t,\bx)} \bu_{\bx}(t,\bx,z)\,{\rm d}z,
\end{align}
giving the evolution for the domain shapes defined in \eqref{eqn:geometricalshapes}. A similar idea is used in boundary integral methods to solve the Stokes problem~\cite{HeltaiArroyoDeSimone-2013-a,ArroyoHeltaiMillan-2012-a}.
The lubrication approximation now provides a simplified energy $E_{\textrm{thin-film}}(h)\equiv\mathcal{E}(h)$ and a simplified operator $K(h)$ such that we directly obtain the evolution equation
\begin{align}
\partial_t h = -K(h){\rm D}\mathcal{E}(h).
\end{align}
For details concerning the construction of the Stokes gradient structure and the lubrication approximation based on energetic arguments we refer to \cite{peschka2018variational}.

The goal of this work is the construction of a rather general class of operators $K(h)$ for the thin-film problem, that also encode the boundary dissipation $\mu_{\partial\omega}$ from \eqref{eqn:newtondiss} systematically into the problem. With the main ideas of the approach encoded in the operator $K$ and the functional $\mathcal{E}$ as defined before, we focus on the construction of efficient numerical schemes for this problem and provide examples showing the necessity to take contact-line dissipation into account. A similar study for the full Navier-Stokes problem was performed by Montefuscolo et al.~\cite{montefuscolo2014high}. However, when one is interested in solving the free boundary problem in the partial-wetting regime {or complete regimes}, then reduced thin-film problems provide quantitative predictions at a significantly lower computational cost.
{
Depending on the magnitude of slip, M\"unch et al.~\cite{munch2005lubrication} showed that there exists a hierarchy of different thin-film models that describes the effective dynamics in these asymptotic regimes. Here, we establish the existence of 
a different hierachy of thin-film models that depends on the magnitude of contact line dissipation.}

First, in section~\ref{sec:thinfilmmodel} we introduce the general gradient flow formalism and define  its ingredients explicitly. Furthermore, we construct the ALE formulation in Eulerian and ALE coordinates and propose a general solution strategy for the time discretisation. In particular, we explore the reduced model that emerges in the strong contact line dissipation limit. In section~\ref{sec:higherorder} we devise higher-order discretisation strategies and discuss the discretisation of line tension to stabilize the motion of the free boundary. Finally, in section~\ref{sec:results} we benchmark the implementations of these free boundary problems and discuss the impact that contact line dissipation has on the observed shapes of sliding droplets and on the instability of a viscous liquid ridge.

\section{Thin-film model and mathematical structures}
\label{sec:thinfilmmodel}
\subsection{General gradient flows}
A gradient flow structure (also known as Onsager structure) is a quadruple 
$(\mathbb{Q},V,\mathcal{E},\Psi^*)$ consisting of a state space $q\in\mathbb{Q}$, a rate space $\partial_t q \simeq \dot{q}\in V$, an energy $\mathcal{E}:\mathbb{Q}\to\mathbb{R}$, and a quadratic dual dissipation $\Psi^*(q,\cdot):V^*\to\R$ defined as 
\begin{align}
  \Psi^*(q,\eta)=\tfrac12\langle \eta,K(q)\eta\rangle,
\end{align}
using the Onsager operator $K(q):V^*\to V$. 
With the dual pairing $\langle\cdot,\cdot\rangle:V^*\times V\to\mathbb{R}$, the usual dissipation $\Psi(q,\cdot):V\to\mathbb{R}$ is defined as the Legendre transform of $\Psi^*$ with respect to the second argument by
\begin{align}
  \Psi(q,\nu)=\sup_{\eta\in V^*}\big(\langle\eta,v\rangle-\Psi^*(q,\eta)\big)=\tfrac12\langle K^{-1}(q)\nu,\nu\rangle,
\end{align}
where $\Psi^*(q,\cdot)$ is convex and $\Psi^*(q,0)\equiv 0$. 
These ingredients define the 
the gradient flow evolution $q:[0,T]\to\mathbb{Q}$ by
\begin{align}
\label{eqn:abstractgradientflowevolution}
  \dot{q}=-K(q){\rm D}\mathcal{E}(q)=-\nabla_{\Psi}\,\mathcal{E}(q),
\end{align}
where the gradient $\nabla_{\Psi}$ is generated by the scalar product $u,v\to \langle K^{-1}(q)u,v\rangle$. While this evolution can be extended to general convex dissipation potentials, here we restrict to the quadratic case defined by $K(q)$. Note that, exploiting the chain rule, these statements directly imply the energy-dissipation balance
\begin{align}
  \tfrac{{\rm d}}{{\rm d}t}\,\mathcal{E}(q(t)) = -2\Psi^*\big(q,{\rm D}\mathcal{E}(q)\big) =-2\Psi(q,\dot{q})\le 0.
\end{align}

In the following, we will show how different choices of  $(\mathbb{Q},V,\mathcal{E},\Psi^*)$ lead to different laws for dynamic contact angles and how they can be discretized naturally. Next, in order to give this abstract definitions a concrete mathematical and physical meaning in terms of partial differential equations, we characterize the state space $\mathbb{Q}$, the energy $\mathcal{E}$, the forces ${\rm D}\mathcal{E}$, and the dual dissipation $\Psi^*$ in more detail.

\subsection{State space, energy, and duality}
Consider for example a system describing the shape of a liquid droplet by an abstract state $q\in\mathbb{Q}$ given by
\begin{align}
  \label{eqn:statespace}
  q=\begin{pmatrix}\omega\subset\mathbb{R}^d\\h:\omega\to\mathbb{R}_+\end{pmatrix}\in\mathbb{Q},
\end{align}
where $h$ denotes the \emph{height} of the liquid droplet over a planar solid surface and $\omega$ is the \emph{support} of the droplet, i.e., outside the domain $\omega\subset\R^d$ liquid is absent and $h\equiv 0$, cf. Fig.~\ref{fig:cangle}. The goal of the manuscript is to derive coupled evolution laws for $q(t)$ and to study corresponding higher-order approximations in \emph{space} and \emph{time}. For the moment we assume $\omega$ is simply connected, $h(x)\ge0$, and $h(x)\equiv 0$ if and only if $x\in\partial\omega$. Furthermore, both $h$ and the boundary $\partial\omega$ are sufficiently smooth.
The volume $\mathit{vol}$ and the area of the support $|\omega|$ of a droplet $q\in\mathbb{Q}$ are given by
\begin{align}
  \mathit{vol}(q)=\int_\omega h\,{\rm d}\omega, \qquad |\omega|(q)=\int_\omega 1{\rm d}\omega. 
\end{align}
In this work we consider conserved dynamics, i.e., with initial data $q_0=q(t=0)$ we have $\mathit{vol}\bigl(q(t)\bigr)\equiv\mathit{vol}(q_0)$.
For some sufficiently smooth 
given $e:\mathbb{R}^d\times\mathbb{R}\times\mathbb{R}^d\to\mathbb{R}$ we then define the energy functional $\mathcal{E}:\mathbb{Q}\to\mathbb{R}$  as
\begin{align}
  \label{eqn:energy}
  \mathcal{E}(q)=\int_\omega e[q] \,{\rm d}\omega, \qquad\text{with}\quad e[q]=e(x,h,\nabla h).
\end{align}
Note that while the energy \eqref{eqn:energy} only depends on the shape of the droplet by \eqref{eqn:statespace}, viscoelastic fluids would additionally require to keep track of elastic deformations in addition to the shape of droplets $q$.
A typical example for the energy of a wetting fluid is $e(x,h,z)=\tfrac\sigma2 |z|^2 + s(x)+ hf(x,h)$, where $\tfrac\sigma2|z|^2+s$ encodes surface energies and $f$ contains bulk energies such as gravity, i.e., $f(\boldsymbol{x},h)=\varrho(g_x\cdot \boldsymbol{x}+\tfrac12g_zh)$ with $\boldsymbol{g}=(g_x,g_z)$ representing the components of gravity in $(\bx,z)$ direction. In the special case where the surface energy is homogeneous $s(\bx)=s\in\R$ and no bulk-energies are present $f\equiv 0$ we can also write the energy functional as $\mathcal{E}=\int \tfrac\sigma2|\nabla h|^2\,{\rm d}x + s|\omega|$. This energy density emerges from \eqref{eqn:newtonenergy} with $s(x)=S=\sigma_\text{SV}-(\sigma_\text{SL}+\sigma)$.
{
\begin{remark} The energy density $e=\frac{\sigma}{2} |z|^2 +s+hf$ is a lubrication approximation of the full energy of the system, where in particular the first term $\frac{\sigma}{2} |z|^2$ is an approximation of the full surface area measure $\sigma(1+|z|^2)^{1/2}$, e.g. \cite{thiele2013gradient}.
\end{remark}
}

Now consider time-dependent curves $q(t)\in \mathbb{Q}$ and the evolution of the energy $\mathcal{E}(q(t))$. Assuming $\dot{h}$ and $h$ are sufficiently regular, the application of Reynold's transport theorem on $\mathcal{E}(q(t))$ gives
\begin{align}
  \nonumber\tfrac{{\rm d}}{{\rm d}t}\,\mathcal{E}(q(t))  &= \int_{\omega} e_z\cdot \nabla\dot{h} + e_h\dot{h}\,{\rm d}\omega + \int_{\partial\omega}e(\dot{\bx}\cdot\bnu)\,{\rm d}\gamma,              
\end{align}
with $\dot{\bx}:\partial\omega\to\mathbb{R}^d$ being the velocity of the boundary which emerges when applying the Reynolds transport theorem and $\bnu:\partial\omega\to\mathbb{R}^d$ is the outer normal vector field. 
In addition to the preceeding construction, a \emph{kinematic condition} needs to be satisfied as an essential constraint.
Therefore, assume $\bx(t)\in\partial\omega(t)$ is any given sufficiently smooth curve on the time-dependent boundary of $\omega(t)$. Due to the virtue of $h(t,\bx(t))\equiv 0$ on $\partial\omega(t)$,  after differentiation with respect to time this becomes
\begin{align}
\label{eqn:kinematiccondition}
  \dot{h}=-\dot{\bx}\cdot\nabla h=\dot{x}|\nabla h|,\qquad \text{on }\partial\omega,
\end{align}
where we use $\dot{x}=\dot{\bx}\cdot\bnu$ for the normal velocity of $\partial\omega$.
While it would be instructive to use Lagrange multipliers to enforce this constraint, due to its simple algebraic form it can be entirely eliminated from the problem of computing $\dot{q}$ by replacing $\dot{x}=\dot{h}/|\nabla h|$ and by modification of the  identification of elements $\eta\in V^*$. 
This allows us to rewrite $\tfrac{{\rm d}}{{\rm d}t}\,\mathcal{E}$ entirely in terms of $\dot{h}$ as 
\begin{align}
\label{eqn:DE}
\tfrac{{\rm d}}{{\rm d}t}\,\mathcal{E}(q(t))  &= \langle {\rm D}\mathcal{E},\dot{h}\rangle=\int_{\omega} e_z\cdot \nabla\dot{h} + e_h\dot{h}\,{\rm d}\omega + \int_{\partial\omega}\frac{e}{|\nabla h|} \dot{h}\,{\rm d}\gamma,
\end{align}
and thereby uniquely defines the driving force ${\rm D}\mathcal{E}$. Next, we define the corresponding dissipation potentials $\Psi,\Psi^*$.

\subsection{Dissipation and weak formulation with kinematic condition}
Assume $h,\omega$ are sufficiently smooth and let $u=(u_\pi\equiv\pi,u_\zeta\equiv\zeta)\in H^1(\omega)\times L^2(\partial\omega)=U$. Define the symmetric, positive bilinear form $a:U\times U\to\R$ by 
\begin{align}
\label{eqn:psiinv}
  a(u,u)=\int_\omega m[q]\,|\nabla \pi|^2\,{\rm d}\omega + \int_{\partial\omega} n[q]\,|\nabla h|^2\,\zeta^2\,{\rm d}\gamma,
\end{align}
and the corresonding invertible operator $A:U\to U^*$ by $a(u_1,u_2)\equiv\langle Au_1,u_2\rangle_U$. Note that both $a=a(q)$ and $A=A(q)$ depend on $q$ through the non-negative coefficient functions $m,n$ and the shape of $\omega$. Furthermore, for $V=H^1(\omega)$ we define the linear operator $L:V\to U^*$ by 
\begin{align}
\label{eqn:Loperator}
  \langle L(v),u \rangle_{U^*\times U} = (v,u_\pi)_{L^2(\omega)} + (v,u_\zeta)_{L^2(\partial\omega)},
\end{align}
with $(u,v)_{L^2(\omega)}=\int_\omega uv\,{\rm d}\omega$ and $(u,v)_{L^2(\partial\omega)}=\int_{\partial\omega} uv\,{\rm d}\gamma$.
This allows us to define the dissipation potential $\Psi:V\to\R$ as
\begin{align}
\Psi(q,\dot{h})=\tfrac{1}{2}\langle L(\dot{h}),A^{-1} L(\dot{h}) \rangle_{U^*\times U}.
\end{align}
Now we are going to define the gradient flow evolution problem $\dot{h}=-K(h){\rm D}\mathcal{E}(h)$ that emerges from the minimisation problem $\dot{h}=\arg\min_v \Psi(q,v)+\langle {\rm D}\mathcal{E},v\rangle_{V^*\times V}$. By differentiation at $\dot{h}$ in the direction of $v$ we get 
$\langle L(v),A^{-1}L(\dot{h})\rangle_{U^*\times U} + \langle {\rm D}\mathcal{E},v\rangle_{V^*\times V}=0$. When we abbreviate $u=-A^{-1}L(\dot{h})\in U$ we can write this as the saddle point problem, where we seek $(u,\dot{h})\in U\times V$ such that
\begin{subequations}
  \label{eqn:derivation-saddle-point}
\begin{align}
  \langle Au,v_u\rangle_{U^*\times U} + \langle L(\dot{h}),v_u\rangle_{U^*\times U}&=0,\\
  \langle L^*(u),v\rangle_{V^*\times V} &= \langle {\rm D}\mathcal{E},v\rangle_{V^*\times V},
\end{align}
\end{subequations}
for all $(v_u,v)\in U\times V$.
Here we used $\langle L^*(u),v\rangle_{V^*\times V}:=\langle L(v),u\rangle_{U^*\times U}$. Here we can use again that $\langle Au,v_u\rangle_{U^*\times U}=a(u,v_u)$ using the definition in \eqref{eqn:psiinv} and the definition of $L$ in \eqref{eqn:Loperator}. Given a solution of this problem exists, then we have the corresponding dual dissipation $\Psi^*(q,\eta)=\tfrac12 a(u,u)$. This allows us to state the final Eulerian weak formulation of the thin-film problem with $u=(\pi,\zeta)$, $v=v_h$, and $v_u=(v_\pi,v_\zeta)$. { Similar saddle-point structures for  fluid-structure interaction problems with inertia and bulk-interface coupling are considered in \cite{peschka2022}.}
\begin{defn}[Weak form of thin-film model]

    The weak form of the thin-film model is to find $(\dot{h},\pi,\zeta)$ such that
    \begin{subequations}
      \label{eqn:weakform}
      \begin{align}
        \label{eqn:weakform1}
          &\int_\omega \pi v_h\,{\rm d}\omega + \int_{\partial\omega} \zeta v_h\,{\rm d}\gamma=\int_{\omega} e_z\cdot\nabla v_h + e_h v_h\,{\rm d}\omega + \int_{\partial\omega}\frac{e}{|\nabla h|}\,v_h\,{\rm d}\gamma, \\
         \label{eqn:weakform2} &\int_\omega \dot{h}v_\pi\,{\rm d}\omega + 
         \int_\omega m[q]\nabla\pi\cdot\nabla v_\pi\,{\rm d}\omega=0,\\
        \label{eqn:weakform3}
        &\int_{\partial\omega}  \dot{h}v_\zeta\,{\rm d}\gamma  
        + \int_{\partial\omega} n[q]\,|\nabla h|^2\,\zeta v_\zeta\,{\rm d}\gamma=0,
      \end{align}
    \end{subequations}
    holds for all test functions $(v_h,v_\pi,v_\zeta)$. 
     In the first equation \eqref{eqn:weakform1} the derivative of the energy is identified with generalized forces ${\rm D}\mathcal{E}\rightsquigarrow u=(\pi,\zeta)$ in the sense $L^*(u)={\rm D}\mathcal{E}$, where  $\zeta:\partial\omega\to\mathbb{R}$ is the force acting on the boundary. Then, in the second equation \eqref{eqn:weakform2} the operator $K$ is applied to the bulk force $\pi$ and in \eqref{eqn:weakform3} it is applied to the boundary force $\zeta$ to obtain the evolution $\dot{h}$ with a dynamic boundary condition.
  \end{defn}
  Note that the expected dependence of the mobilities on the state is given by $m[q]\equiv m(h)$ and $n[q]\equiv n(|\nabla h|)$, e.g., the dependence of the degenerate bulk mobility $m$ on the height $h$ is $m(h)=\tfrac{1}{3\mu}h^3+bh^2$ with $b$ the Navier-slip length parameter. For heterogenous or switchable substrates one can also have $m\equiv m(t,\bx,h)$ and $n\equiv n(t,\bx,|\nabla h|)$. Since we do not supply an external force to the gradient evolution, this dependence of $m,n$ on time and space does not violate the energetic structure, i.e., $\tfrac{\rm d}{{\rm d}t}\mathcal{E}(q(t))=-2\Psi^*\bigl(t,q,{\rm D}\mathcal{E}(q)\bigr)\equiv -2\Psi(t,q,\dot{h})\le 0$.

\subsection{Weak formulation in Arbitrary Lagrangian Eulerian form}

The formulation of the thin-film problem using the state variables $q$ as introduced in \eqref{eqn:statespace} is very elegant but computationally impractical. In the following we introduce the Arbitrary Lagrangian Eulerian form of the evolution equations \eqref{eqn:weakform}. This formulation is based on a reference domain $\bar\omega\subset\R^d$ that does not depend on time, e.g., often one uses the initial support $\bar\omega=\omega(t=0)$.

Then one introduces the ALE state $\bar{q}$ using diffeomorphisms and height functions defined on the reference domain as
\begin{subequations}
\label{eqn:ale_statespace}
\begin{align}
  \bar{q}=\begin{pmatrix}\bar{\psi}:\bar\omega\to\mathbb{R}^d\\\bar{h}:\bar{\omega}\to\mathbb{R}_+\end{pmatrix}\in\bar{\mathbb{Q}}.
\end{align}
This state variable relates to the original state variable \eqref{eqn:statespace} by the conditions
\begin{align}
\label{eqn:ale_domain}&\omega=\bar\psi(\bar\omega)\equiv\{\bx\in\R^d:\exists\bar{\bx}\in\bar\omega\text{ s.t. } \bx=\bar\psi(\bar{\bx})\},\\
\label{eqn:ale_height}&h(\bar\psi(\bar{\bx}))=\bar{h}(\bar{\bx}).
\end{align}
\end{subequations}
The description of the state $q$ by a state $\bar{q}$ is by no means unique but there exists a class of equivalent descriptions in the following sense:
\begin{remark}[Equivalent ALE descriptions]
Let $\bar{q}=(\bar{\psi},\bar{h})$ an ALE state as defined \eqref{eqn:ale_statespace} in describing a state $q=(\omega,h)$ introduced in \eqref{eqn:statespace} and consider any diffeomorphism $\hat{\psi}:\tilde{\omega}\to\bar\omega$ using a new set $\tilde\omega\subset\R^d$. Then 
\begin{align*}
  \tilde{q}=\begin{pmatrix}\tilde{\psi}:\tilde\omega\to\mathbb{R}^d\\\tilde{h}:\tilde{\omega}\to\mathbb{R}_+\end{pmatrix},
\end{align*}
is an equivalent ALE state, when we define the new diffeomorphism by $\tilde\psi=\bar\psi\circ\hat\psi$ and correspondingly $h(\tilde{\psi}(\tilde{x}))=\tilde{h}(\tilde{x})$. Even if $\tilde\omega=\bar\omega$ this equivalence class contains diffeomorphisms that map $\bar\omega$ to itself and thereby preserve the shape of the domain. The goal of the Arbitratry Langrangian Eulerian framework to select a trajectory of ALE descriptions within this equivalence class, such that the evolution equations for $q:[0,T]\to\mathbb{Q}$ are replaced by suitable evolution equations for $\bar{q}:[0,T]\to\bar{\mathbb{Q}}$ satisfying \eqref{eqn:ale_statespace}.
\end{remark}

We introduce the ALE deformation gradient $\bar{F}:\bar\omega\to\R^{d\times d}$, the Jacobian $\bar{J}:\bar\omega\to\R$, and the area $\bar{A}:\partial\bar\omega\to\R$
\begin{align}
  \bar{F}(\bar{x})=\bar{\nabla}\bar\psi(\bar{\bx}), \qquad \bar{J}(\bar{\bx})=\det\bar{F}(\bar{\bx}), \qquad \bar{A}(\bar{\bx})=\bar{J}\|\bar{F}^{-\top}\bar\nu\|,
\end{align}
using the normal vector field $\bar\nu:\partial\bar\omega\to\R^d$.
By computing the Eulerian time-derivative of $h$ in \eqref{eqn:ale_height} we obtain the convective derivative 
\begin{align}
\label{eqn:convective_derivative}
\dot{\bar{h}}(t,\bar{\bx})=\dot{h}(t,\bx)+\dot{\bar{\psi}}(t,\bar{\bx})\cdot\nabla h(t,\bx),
\end{align}
with $\bx=\bar{\psi}(t,\bar{\bx})$.

\begin{defn}[Weak ALE form of thin-film model]
\label{eqn:weakALEformulation}
In the following we present two variants of the weak ALE form of the thin-film model, i.e., where the time derivative of the mapping $\dot{\bar{\psi}}$ is coupled i) or decoupled ii) from the weak formulation. While i) and ii) are equivalent on the continuous level, the discretisation of either variant comes with certain advantages and disadvantages.

\noindent
\begin{enumerate}[i)]
\item\textit{Coupled formulation:} The weak form of the thin-film model is to find the quadruple of functions $(\dot{\bar{h}}:\bar{\omega}\to\R,\bar{\pi}:\bar\omega\to\R,\bar{\zeta}:\partial\bar\omega\to\R,\dot{\bar{\psi}}:\bar\omega\to\R^d)$ such that
  \begin{subequations}
    \label{eqn:ale_weakform}
    \begin{align}
      \label{eqn:ale_weakform1}
        &\int_{\bar\omega} \bar{\pi} v_h\,{\bar{J}\rm d}{\bar\omega} + \int_{\partial{\bar\omega}} \bar{\zeta} v_h\,\bar{A}{\rm d}\bar\gamma=\int_{{\bar\omega}} e_z\cdot\bar{F}^{-\top}\bar{\nabla}v_h + e_h v_h\,{\bar{J}\rm d}{\bar\omega} + \int_{\partial{\bar\omega}}\frac{e}{|\bar{F}^{-\top}\bar{\nabla} h|}\,v_h\,\bar{A}{\rm d}\bar\gamma, \\
       \label{eqn:ale_weakform2} &\int_{\bar\omega} \dot{h}\bar{v}_\pi\,{\bar{J}\rm d}{\bar\omega} + 
       \int_{\bar\omega} m[q]\bar{F}^{-\top}\bar{\nabla}\bar{\pi}\cdot\bar{F}^{-\top}\bar{\nabla}\bar{v}_\pi\,{\bar{J}\rm d}{\bar\omega}=0,\\
      \label{eqn:ale_weakform3}
      &\int_{\partial{\bar\omega}}  \dot{h}\bar{v}_\zeta\,{\bar{A}\rm d}\bar\gamma  
      + \int_{\partial{\bar\omega}} n[q]\,|\bar{F}^{-\top}\bar{\nabla}\bar{h}|^2\,\bar{\zeta} \bar{v}_\zeta\,{\bar{A}\rm d}\bar\gamma=0,
    \end{align}
  
  for all test functions $(\bar{v}_h,\bar{v}_\pi,\bar{v}_\zeta,\bar{v}_\psi)$ but is still formulated in terms of $\dot{h}$ and $v_h$. Using \eqref{eqn:convective_derivative}, one replaces $\dot{h}=\dot{\bar{h}}-\dot{\bar{\psi}}\cdot\bar{F}^{-\top}\bar{\nabla}\bar{h}$ and correspondingly $v_h=\bar{v}_h-\bar{v}_\psi\cdot\bar{F}^{-\top}\bar{\nabla}\bar{h}$ in \eqref{eqn:ale_weakform} to obtain the weak formulation corresponding to \eqref{eqn:weakform}. Observe that while the Eulerian time-derivative $\dot{h}$ has natural boundary conditions, the ALE time derivative needs to vanish on the boundary, i.e., $\dot{\bar{h}}=0$ and $\bar{v}_h=0$ on $\partial\bar\omega$. This is taken into account in the weak formulation as an essential boundary condition and corresponding terms do not appear in boundary integrals upon replacement of $\dot{h}$ and $v_h$. In order to complete this description we need to provide the weak formulation for the rate of the displacement $\dot{\bar{\psi}}$. One possible choice is
  \begin{align}
  \label{eqn:ale_weakform4}
    \int_{\bar\omega} \mathbb{D}\dot{\bar{\psi}}:\mathbb{D}\bar{v}_\psi\,{\bar{J}\rm d}\bar\omega = 0,
  \end{align}
  using the symmetric gradient $\mathbb{D}\bar{v} = \tfrac12\big(\bar\nabla\bar{v}\bar{F}^{-1}+\bar{F}^{-\top}(\bar\nabla\bar{v})^\top\big)$. The boundary conditions for $\dot{\bar{\psi}}$ emerge from the boundary terms in \eqref{eqn:ale_weakform}. 
  The remaining ambiguity of the tangential rate at the boundary can be used to select from the equivalence class of ALE description a representation that ensures the quality of the transformation upon discretisation. However, often the tangential component is simply set to zero.
\end{subequations}
\item\textit{Decoupled formulation:} The basis of the decoupled formulation is the observation that equations \eqref{eqn:ale_weakform} 
only depend on $\dot{\bar{\psi}}$ upon the replacement of $\dot{h},v_h$ using \eqref{eqn:convective_derivative} as explained in the coupled formulation. While these equations depend nonlinearly on $\psi$, in this form they do not depend on the rate $\dot{\psi}$. Thus we can devise the following solution strategy in three steps.\\[-1em]

\noindent
\textbf{Step 1:} In the first step of the decoupled formulation  one seeks $(\dot{h}:\bar{\omega}\to\R,\bar{\pi}:\bar\omega\to\R,\bar{\zeta}:\partial\bar\omega\to\R)$ such that \eqref{eqn:ale_weakform}
are valid for all test functions $(v_h,\bar{v}_\pi,\bar{v}_\zeta)$. Note that here $\dot{h}$ refers to the Eulerian time-derivative evaluated on the reference domain $\bar\omega$, where for convenience we use the same symbol as for the original Eulerian time-derivative.\\[-0.6em]

\noindent
\textbf{Step 2:} The second step is concerned with the reconstruction of the deformation in $\dot{\bar{q}}=(\dot{\bar{\psi}},\dot{\bar{h}})$ from the $\dot{h}$ computed in the first step. Therefore, first notice that on the boundary $\partial\bar\omega$ we have $0=\dot{\bar{h}}=\dot{h}+\dot{\bar{\psi}}\cdot\bar{F}^{-\top}\bar{\nabla}\bar{h}$. Thus, using the previously computed $\dot{h}$ we obtain the (Eulerian) normal component of $\dot{\psi}$ on the boundary via 
\begin{subequations}
\label{eqn:decoupled_formulation}
\begin{align}
\label{eqn:normalcomponents}
\nu\cdot\dot{\bar{\psi}}=\frac{\dot{h}}{|\bar{F}^{-\top}\bar{\nabla}\bar{h}|}, \qquad \nu=-\frac{\bar{F}^{-\top}\bar{\nabla}\bar{h}}{|\bar{F}^{-\top}\bar{\nabla}\bar{h}|}.
\end{align}
The tangential component can be chosen arbitrarily, e.g., $\dot{\bar{\psi}}\cdot t=0$ is a common but not necessarily optimal choice. Using the boundary values determined from $\dot{h}$ we can now solve \eqref{eqn:ale_weakform4} to determine $\dot{\bar{\psi}}$, specifically we now solve for $(\dot{\bar{\psi}}:\bar\omega\to\R^d,\bar\lambda:\partial\bar\omega\to\R^d)$ such that
\begin{align}
&\int_{\bar\omega} \mathbb{D}\dot{\bar{\psi}}:\mathbb{D}\bar{v}_\psi\,\bar{J} {\rm d}\bar\omega + \int_{\partial\bar\omega}\bar{\lambda}\cdot\bar{v}_\psi\,{\rm d}\bar\gamma= 0,\\
&\int_{\partial\bar\omega} \dot{\bar{\psi}}\cdot\bar{v}_\lambda \,{\rm d}\bar\gamma = \int_{\partial\bar\omega} -\bar{v}_\lambda\cdot\left(\frac{\bar{F}^{-\top}\bar{\nabla}\bar{h}}{|\bar{F}^{-\top}\bar{\nabla}\bar{h}|^2}\dot{h}+c_1 t\right)\,{\rm d\bar\gamma},
\end{align}
for all $(\bar{v}_\psi,\bar{v}_\lambda)$. Here $\bar\lambda$ is a Lagrange multiplier for the constraint enforcing the Dirichlet boundary conditions and $c_1$ contains a possible (arbitrary) correction of the tangential component of $\dot{\bar{\psi}}$.\\[-0.6em]

\noindent
\textbf{Step 3:} In the third and final step we compute $\dot{\bar{h}}$ using the Lagrangian variant of \eqref{eqn:convective_derivative}, i.e.,
\begin{align}
\dot{\bar{h}}(t,\bar{\bx})=\dot{h}(t,\bar{\bx})+\dot{\bar{\psi}}(t,\bar{\bx})\cdot\bar{F}^{-\top}\bar{\nabla}\bar{h}(t,\bar{\bx}),
\end{align}
\end{subequations}
using $\dot{h}$ computed in step 1 and using $\dot{\bar{\psi}}$ computed in step 2. The resulting $\dot{\bar{h}}$ automatically satisfies the desired boundary condition $\dot{\bar{h}}=0$ on $\partial\bar\omega$ by construction.
\end{enumerate}
{
The coupled  and decoupled formulation \eqref{eqn:ale_weakform} and \eqref{eqn:decoupled_formulation} of the continuous thin-film problem with dynamic contact angle solve the same problem since both emerge from the same weak formulation \eqref{eqn:weakform} upon a change of variables from Eulerian time derivative $\dot{h}$ to an ALE time-derivative $\dot{\bar{h}}$. However, once one introduces a time-discretization both formulations can potentially differ due to a different choice of semi-implicit discretization. In such a case, the coupled formulation might be better suited if one desires to implement a fully implicit time-discretisation. 
}
Alternatively, the decoupled formulation will in general be easier to implement because solving for $\dot{h}$ in step 1 and the reconstruction of $\dot{\bar{q}}=(\dot{\bar{\psi}},\dot{\bar{h}})=\mathcal{F}(\bar{q})$ are decoupled and lead to smaller systems of linear equations. Then the evolution of the droplet is computed by solving
\begin{align}
  \dot{\bar{q}}=\mathcal{F}(\bar{q}),
\end{align}
using a suitable discretisation. The original state can then be reconstructed via \eqref{eqn:ale_statespace}.
\end{defn}

\subsection{Weak and strong contact line dissipation limits}
\paragraph{Thin-film model with dynamic contact line}
Before we have stated the thin-film models in their variational gradient flow
form. This is extremely useful if one is interested in the discretisation of
these equations using finite elements. Also, the transition from the thin-film
PDE to the weak formulation is not at all obvious compared to 
the reverse transition from the weak formulation to the PDE. Therefore, using the weak formulation in Eulerian coordinates provided
in \eqref{eqn:weakform}, we assume the the solution is sufficiently regular to
employ integration by parts. Starting with \eqref{eqn:weakform2} and for given
initial data $q_0=(\omega(t=0),h(t=0))$ we obtain
\begin{subequations}
\label{eqn:strong_thinfilm}
\begin{align}
\label{eqn:thinfilm_strongform}
  \dot{h}-\nabla\cdot (m[q]\nabla\pi)=0,\qquad&\text{in }\omega,\\
  \label{eqn:thinfilm_noflux}
  m[q]\nu\cdot\nabla\pi=0,\qquad&\text{on }\partial\omega,
\end{align}
where the second natural boundary condition ensures the conservation of volume $\tfrac{\rm d}{{\rm d}t}\int_\omega h\,{\rm d}x=0$. Using, for example, the energy density $e[q]=\tfrac\sigma2|\nabla h|^2+s+hf(x,h)$ we obtain from \eqref{eqn:weakform1} the generalized forces/pressures
\begin{align}
\label{eqn:force_pi}
\pi=-\sigma\Delta h + \partial_h (hf),\qquad&\text{in }\omega,\\
\label{eqn:force_zeta}
\zeta=|\nabla h|^{-1}\big(-\tfrac{\sigma}{2}|\nabla h|^2+s\big),\qquad&\text{on }\partial\omega.
\end{align}
Finally, using \eqref{eqn:weakform3} we directly obtain the dynamic boundary condition (dynamic contact line)
\begin{align}
\label{eqn:thinfilm_dynamicangle}
\dot{h}+n[q]|\nabla h|^2\zeta=0,\qquad&\text{on }\partial\omega.
\end{align}
These equations are supplemented with the kinematic condition for the motion of the domain, i.e.,
\begin{align}
\label{eqn:thinfilm_kinematic}
\dot{h}+\dot{\mathbf{x}}\cdot\nabla h=0,\qquad&\text{on }\partial\omega.
\end{align}
\end{subequations}
The equations in \eqref{eqn:strong_thinfilm} are the full thin-film model with dynamic contact angle as a dynamic boundary condition and a general driving force encoded in the energy $e$. If desired, one can insert the forces $\pi,\zeta$ from (\ref{eqn:force_pi},\ref{eqn:force_zeta}) into the partial differential equation \eqref{eqn:thinfilm_strongform} and \eqref{eqn:thinfilm_dynamicangle}. Alternatively, one can also represent the dynamic contact angle \eqref{eqn:thinfilm_dynamicangle} as an equation for the normal component of the velocity $\dot{x}=\dot{\mathbf{x}}\cdot\nu=-\dot{\mathbf{x}}\cdot\nabla h/|\nabla h|$, i.e.,
\begin{align}
  \label{eqn:thinfilm_dynamicangle1}
  \dot{x}=-n[q]\,|\nabla h|\,\zeta =-n[q]\,\big(-\tfrac{\sigma}{2}|\nabla h|^2+s\big).
\end{align}
Often one has a mobility $m[q]={m}_0|h|^\alpha$ with ${m}_0>0$ which supports a moving contact line for $0\le \alpha <3$ but produces a logarithmic singularity for $\alpha=3$. Similarly, for the dynamic contact angle we have $n[q]=n_0 |\nabla h|^\beta$ for $n_0>0$. The equilibrium contact slope (angle) implied by \eqref{eqn:thinfilm_dynamicangle1} is $|\nabla h|=(2s/\sigma)^{1/2}$ for $s\ge 0$. 

{
\begin{remark}[Existing work on dynamic contact angles]
Dynamic contact angles for thin films have been known for a long time \cite{cox1986dynamics,voinov1976hydrodynamics,hocking1983spreading} and are attributed to localized dissipation at the contact line \cite{de2013capillarity}. Several works consider hydrodynamic models and microscopic origins of dynamic contact angle, e.g. \cite{snoeijer2013moving,bonn2009wetting,shikhmurzaev1997moving} and 
can be directly applied to thin film models in a small slope approximation $\sin\vartheta\approx|\nabla h|\ll 1$. 

Only a few studies go beyond the formulation of laws for dynamic contact angle and treat the moving contact line as a sharp-interface. Limat et al. \cite{snoeijer2007cornered} consider sliding droplets with dynamic contact angle using a formal expansions around a singular corner tip. Works by Kn\"upfer, Gnann, Otto \& Giacomelli \cite{knupfer2015well,giacomelli2014well} perform the mathematical analysis of models for spreading or partial wetting and in particular Chiricotto \& Giacomelli \cite{chiricotto2017weak} recently proposed a model for a dynamic contact angle based on the hydrodynamic work by Ren \& E \cite{ren2011derivation}. A related model was derived by Doi et al. \cite{xu2016variational} based on Onsager's principle and studied numerically used a robust approximation. Sharp-interace limits are derived using formal asymptotics by King et al. \cite{king2009linear}. Similar fourth-order problems based on Wasserstein gradient flows but without dynamic contact angles have been considered using Lagrangian schemes in \cite{matthes2017convergent}.
However, the authors of this work are not aware of any work that presents a weak formulation for thin films with dynamic contact angle based on gradient flows or implements a full transient thin-film model with dynamic contact angles.
\end{remark}
}

{ In the following we discuss the limits of strong and of weak contact line dissipation, which in terms of mobilities correspond the the limits of vanishing contact line mobility $n\to 0$ and of infinite contact line  mobility $n\to\infty$, respectively.}
We present the resulting system of partial differential equation and, as before, the corresponding weak formulation that is suitable for a finite element discretisation.

\paragraph{Weak contact line dissipation limit $n\to\infty$}
Let us first consider the limit of weak contact line dissipation or large contact line mobility $n\to\infty$. This leads to \eqref{eqn:force_zeta} and \eqref{eqn:thinfilm_dynamicangle} together, resulting in a boundary condition where the contact angle is in equilibrium, i.e., $|\nabla h|=(2s/\sigma)^{1/2}$.

For given initial data $q_0=(\omega(t=0),h(t=0))$ one seeks solutions $h,\pi$ of
\begin{subequations}
  \label{eqn:weakfriction_thinfilm}
  \begin{align}
    \dot{h}-\nabla\cdot (m[q]\nabla\pi)=0,\qquad&\text{in }\omega,\\
    \pi=-\sigma\Delta h + \partial_h (hf),\qquad&\text{in }\omega,\\
    m[q]\nu\cdot\nabla\pi=0,\qquad&\text{on }\partial\omega,\\
    |\nabla h|=(2s/\sigma)^{1/2},\qquad&\text{on }\partial\omega,
  \end{align}
  together with the kinematic condition for the evolution of the domain $\omega(t)$
  \begin{align}
    \dot{h}+\dot{\mathbf{x}}\cdot\nabla h=0,\qquad&\text{on }\partial\omega.
  \end{align}
\end{subequations}
For the system to be consistent, we require that the initial conditions already satisfy the equilibrium contact angle condition $|\nabla h(t=0)|=(2s/\sigma)^{1/2}$. In the corresponding weak formulation of the weak contact line dissipation limit one just needs to adjust the identification of boundary forces. Then, the forces at the contact line are always in equilibrium and in the weak formulation we seek $(\dot{h}:\omega\to\R,\pi:\omega\to\R)$ such that
      \begin{subequations}
        \label{eqn:oldweakform}
        \begin{align}
          \label{eqn:oldweakform1}
            &\int_\omega \pi v_h\,{\rm d}\omega =\int_{\omega} e_z\cdot\nabla v_h + e_h v_h\,{\rm d}\omega + \int_{\partial\omega}\frac{e}{|\nabla h|}\,v_h\,{\rm d}\gamma=\langle{\rm D}\mathcal{E}(q),v_h\rangle, \\
           \label{eqn:oldweakform2} &\int_\omega \dot{h}v_\pi\,{\rm d}\omega +
           \int_\omega m[q]\nabla\pi\cdot\nabla v_\pi\,{\rm d}\omega=0,
        \end{align}
      \end{subequations}
holds true for all test functions $(v_h:\omega\to\R,v_\pi:\omega\to\R)$. This weak formulation is equivalent to the variational formulation of the free boundary thin-film problem presented in \cite{peschka2015thin}. 

\paragraph{Strong contact line dissipation limit $n\to 0$}
Let us finally consider the slightly more subtle limit of strong contact line dissipation or vanishing contact line mobility $n\to 0$, which we perform by setting the mobility to
\begin{align}
n[q]=\varepsilon \hat{n}[q],
\end{align}
and then let $\varepsilon\to 0$. For given initial data $q_0=(\omega_0,h_0)$ with volume $\mathit{vol}(q_0)$ consider a solution $h_\varepsilon$ of \eqref{eqn:strong_thinfilm}. 
What happens for $\varepsilon\to 0$ is that the contact line freezes and from \eqref{eqn:thinfilm_dynamicangle} and \eqref{eqn:thinfilm_kinematic} we deduce $\dot{h}_\varepsilon \to 0$ and $\dot{x}_\varepsilon \to 0$ on $\partial\omega$ for any fixed time $T$. 
Similarly, if we let $T\to\infty$ but $T = o(\varepsilon^{-1})$ as $\varepsilon\to0$, then $h_\varepsilon(T,\boldsymbol{x})$ approaches a stationary solution $\hat{h}$ which minimizes the energy for any fixed given $\hat\omega$, i.e., $h_\varepsilon(T,\boldsymbol{x})\to \hat{h}_{\hat{\omega}}(\boldsymbol{x})$ with
\begin{subequations}
\label{eqn:strongfriction}
\begin{align}
\label{eqn:strongfriction_energy}
  \hat{h}_{\hat{\omega}} = \argmin_{h,\mathit{vol}(q)=\mathit{vol}(q_0)} \mathcal{E}\big(q=(h,\hat{\omega})\big),
\end{align}
provided that such a solution exists and that it is meaningful, which we comment on later. If we rescale the time $\hat{t}=\varepsilon t$ such that $\dot{\hat{x}}=\varepsilon^{-1}\dot{x}$, then the domain evolves again according to
\begin{align}
\label{eqn:strongfriction_dynamic}
\dot{\hat{x}}=-\hat{n}[\hat{q}]\,\big(-\tfrac\sigma2|\nabla\hat{h}_{\hat\omega}|^2+s\big),
\end{align}
\end{subequations}
for the stationary profile $\hat{h}$ from the minimisation \eqref{eqn:strongfriction_energy}. 
For any given initial data $\hat{q}_0=(\hat{\omega}_0,\hat{h}_{0})$, the system \eqref{eqn:strongfriction} is a gradient flow evolution, where $\hat{h}$ is already determined by the initial volume $\mathit{vol}({q}_0)$ and $\hat{\omega}$ through \eqref{eqn:strongfriction_energy}, i.e., $\hat{h}(t)\equiv\hat{h}_{\hat{\omega}(t)}$. The weak formulation for strong contact line dissipation limit is to find $(\hat{\pi}\in\R,\hat{h}:\omega\to\R)$ such that
%
%
%
\begin{subequations}
  \label{eqn:strongfriction_weakform}
  \begin{align}
    &\int_\omega \hat{\pi} v_h\,{\rm d}\omega =\int_{\omega} e_z\cdot\nabla v_h + e_h v_h\,{\rm d}\omega, \\
   &\int_\omega \hat{h}_{\hat\omega}\,{\rm d}\omega =\mathit{vol}(\hat{q}_0),
\end{align}
\end{subequations}
with homogeneous Dirichlet boundary conditions for $\hat{h}$ and the test functions $v_h$ and the Lagrange multiplier $\hat{\pi}\in\R$. From the corresponding solution we determine the contact line velocity using \eqref{eqn:strongfriction_dynamic}. This equation is simpler to solve than \eqref{eqn:weakform} since the bulk force balance and boundary force balance decouple and are now contained in \eqref{eqn:strongfriction_weakform} and \eqref{eqn:strongfriction_dynamic}.

This and similar types of quasi-static droplet evolution have been previously investigated with repect to aspects of modeling, numerics, and mathematical well-posedness, e.g., cf. \cite{greenspan1978motion,glasner2005boundary,grunewald2011variational}. While model hierarchies for quasi-static evolution have been discussed for the Stokes free boundary problem, e.g., cf. \cite{shikhmurzaev1997spreading}, similar hiearchies and in particular the numerical treatment of thin-film free boundary problems with dynamic contact angles are poorly investigated in literature. Despite the lack of a rigorous proof or a formal asymptotic argument, the variational models presented here are the natural candidates for a strong or weak contact line dissipation limit. The following remark emphasizes the need for such a hierarchy.

\begin{remark}[Feasibility of $\hat{h}_{\hat{\omega}}$] 
  \label{rem:feasibility}
    In general not all states $\hat{q}=(\hat{\omega},\hat{h}_{\hat{\omega}})$ that result from the minimisation in \eqref{eqn:strongfriction_energy} are feasible, i.e., it is not guaranteed that minimizers satisfy $\hat{h}_{\hat{\omega}}\ge 0$ inside $\hat{\omega}$. However, violation of this condition renders such a solution meaningless and one would have to fall back to the full problem \eqref{eqn:strong_thinfilm} and expect topological transitions on a shorter time scale. 
    { For a concrete counterexample on a unit disc $\hat{\omega}=B_1(0)$ see Fig.~\ref{fig:counter_example}.}

\begin{figure}[H]
      \centering
      \includegraphics[width=\textwidth]{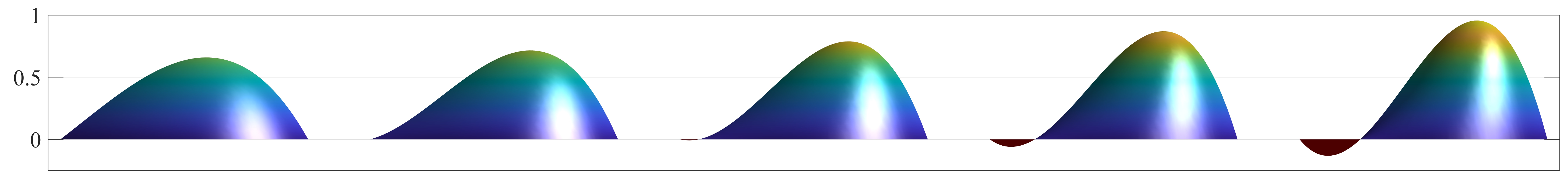}
      \caption{Side view of stationary droplet shapes $\hat{h}_{\hat{\omega}}$ on a unit disc $\hat{\omega}=B_1(0)$ with $\mathit{vol}=1$, $\sigma=1$, $\epsilon=0.01$, $g_z=0$ and gravity in $x$-direction increasing from (left) to (right) $g_x=\{2,4,6,8,10\}$ showing that positivity of solutions  is  is violated for $g_x\ge 6$.}  
      \label{fig:counter_example}
\end{figure}
  
{
    Topological transitions are a general issue for virtually all sharp-interface models for free boundary problems. For related mean curvature flows ``surgery'' as introduced by Huisken \& Sinestrari \cite{huisken2009mean} is a candidate for a well-defined procedure to handle topological transitions for certain thin-film problems, which we leave for future work.}
\end{remark}

{
\subsection{Related thin-film models}
While the classical thin-film model for viscous fluids \cite{oron1997long} has properties such as non-negativity or finite propagation velocity \cite{bernis1996finite}, there is a much larger class of related models that possess a similar gradient structure and thus require similar numerical methods. In the following, we will discuss some similarities and differences in more detail.
A more general class of nonlinear parabolic PDEs is given by 
\begin{subequations}
\begin{align}
\partial_t h - \nabla\cdot(m[h] \nabla \mu) + M[h]\mu=0, \qquad \mu=\frac{\delta E[h]}{\delta h}, 
\end{align}
where $m,M$ are positive mobilities that depend on the state $h$, i.e.,
\begin{align}
  m[h]=m(h,\nabla h)\ge 0, \qquad M[h]=M(h,\nabla h)\ge 0,
\end{align}
and the energy is given by some density
\begin{align}
\mathcal{E}(h)=\int E[h]\,{\rm d}x, \qquad E[h]=E(h,\nabla h,\Delta h),
\end{align}
for some $E:\R\times\R^d\times \R\to\R$ 
such that $\mu = \partial_h E - \nabla\cdot (\partial_{\nabla h} E) + \Delta(\partial_{\Delta h}E)$. For $m>0$ and $M=0$ the dynamics is a $H^{-1}$ gradient flow and for $m=0$ and $M>0$ the dynamics is a $L^2$ gradient flow. Otherwise the dynamics is of mixed type. Smooth solutions satisfy the energy balance
\begin{align}
  \frac{\rm d}{{\rm d}t}\mathcal{E}(h(t)) = -\int m |\nabla\mu|^2 + M\mu^2\,{\rm d}x \le 0.
\end{align}
\end{subequations}

With $E[h]=E(h,z,d)$, then the following classes of models are considered in literature as special cases:
\begin{enumerate}[i)]
\item Cahn-Hilliard models: $E[h]=\tfrac\varepsilon2|z|^2 + W_0(h)$, $m[h]=m_0(h)$, $M\equiv 0$,

\item Allen-Cahn models:
$E[h]=\tfrac\varepsilon2|z|^2 + W_0(h)$, $m\equiv 0$, $M[h]=M_0(h)$,

\item reaction-diffusion: $E[h]=h(\log h - 1)$, $m[h]=m_0(h)$, $M[h]=M_0(h)$,

\item epitaxy: $E[h] = \tfrac\delta 2|d|^2 + W_1(z)$, $m\equiv 0$, $M[h]=M_1(z)$.
\end{enumerate}

Cahn-Hilliard-type models i) with nonconvex $W_0$ are typically used to describe diffusion and phase-separation or become classical thin-film models when $W_0$ is an intermolecular potential or wetting energy. For thin-films usually $m_0(h)=\hat{m}_0 |h|^\alpha$, where the exponent $\alpha\ge 0$ of the degeneracy determines the qualitative behavior of the support, which is the crucial property for this work and makes the numerical treatment challenging, whereas if $h>0$ is ensured, usually standard methods for parabolic PDEs are available. 
Allen-Cahn-type models ii) with nonconvex $W_0$ are typically used to describe phase-transitions. For a convex energy iii) with general mobilities $m_0$ and $M_0$ usually reaction-diffusion systems are described. 
Gradient-dependent mobilities and energies with dependency on second-order derivatives and with Ehrlich–Schwoebel terms are used in epitaxy models with slope selection $W_1(z)=\tfrac14 (|z|^2-1)^2$ or without slope selection $W_1(z)=-\tfrac12 \ln(1+|z|^2)$.
All above models share the gradient structure and therefore energy-stable discretizations are highly desireable, in particular if the long-time or coarsening dynamics is of interest. Research has focussed on time-discretizations, e.g., using convex-concave splitting, and on energy-stability \cite{wang2010unconditionally,zhang2021structure}, or positivity preserving \cite{grun2000nonnegativity,zhornitskaya1999positivity} schemes for global models \cite{diez2000global}. Extensions for anisotropic energies have been considered in \cite{korzec2013anisotropic,wise2005quantum}.
}

\section{Higher-order discretisation of moving domains}
\label{sec:higherorder}
In the following we will present a space and time-discretisation of the thin-film problem with dynamic contact angle. For simplicity we will focus on higher-order discretisations of the decoupled weak ALE formulation introduced in Definition~\ref{eqn:weakALEformulation} and follow the presented 3-step strategy in order to construct a single first-order time step evolution step of the state $\bar{q}^{n-1}=\bar{q}(t)$ 
to $\bar{q}^{n}=\bar{q}(t+\tau)$. 
{
Note that we use $\dot{h}$ for Eulerian time derivatives. Since solutions $h$ at time $t$ and $t+\tau$ are defined on different domains, this expression cannot be approximated by a difference quotient. Therefore, we stick to the unconventional dot-notation and  introduce difference quotients only later in the ALE frame.
}
\subsection{Single decoupled step}
\paragraph*{Step 1} In the first step, for a given $\bar{q}^{n-1}\equiv(\bar{\psi},\bar{h})$ and corresponding explicitly given $\bar{F}=\bar\nabla\bar\psi$, $\bar{J}=\det\bar{F}$ and $\bar{A}=\bar{J}\|\bar{F}^{-\top}\bar{\nu}\|$ we seek 
$\dot{h},\bar{\pi}:\bar{\omega}\to\mathbb{R}$ { and $\bar{\zeta}:\partial\bar{\omega}\to\mathbb{R}$ such that}
  \begin{align*}
      &\int_{\bar\omega} \bar{\pi} v_h\,{\bar{J}\rm d}{\bar\omega} + \int_{\partial{\bar\omega}} \bar{\zeta} v_h\,\bar{A}{\rm d}\bar\gamma=\int_{{\bar\omega}} e^n_z\cdot\bar{F}^{-\top}\bar{\nabla}v_h + e_h v_h\,{\bar{J}\rm d}{\bar\omega} + \int_{\partial{\bar\omega}}\frac{e}{|\bar{F}^{-\top}\bar{\nabla} \bar{h}|}\,v_h\,\bar{A}{\rm d}\bar\gamma, \\
     &\int_{\bar\omega} \dot{h}\bar{v}_\pi\,{\bar{J}\rm d}{\bar\omega} + 
     \int_{\bar\omega} m[\bar{q}]\bar{F}^{-\top}\bar{\nabla}\bar{\pi}\cdot\bar{F}^{-\top}\bar{\nabla}\bar{v}_\pi\,{\bar{J}\rm d}{\bar\omega}=0,\\
    &\int_{\partial{\bar\omega}}  \dot{h}\bar{v}_\zeta\,{\bar{A}\rm d}\bar\gamma  
    + \int_{\partial{\bar\omega}} n[\bar{q}]\,|\bar{F}^{-\top}\bar{\nabla}\bar{h}|^2\,\bar{\zeta} \bar{v}_\zeta\,{\bar{A}\rm d}\bar\gamma=0,
  \end{align*}
  { holds for all $(\nu_h,\bar{\nu}_\pi,\bar{\nu}_\zeta)$. In a later time-discretization, } the highest-order term $e_z$ is treated semi-implicitly in $\dot{h}$ by setting $e^n_z=\partial_z e(\bar{\psi},\bar{h},\bar{F}^{-\top}\bar{\nabla}(h+\tau\dot{h}))$ but the lower-order terms are treated explicitly by setting $e_h=\partial_h e(\bar{\psi},\bar{h},\bar{F}^{-\top}\bar{\nabla}h)$, $e=e(\bar{\psi},\bar{h},\bar{F}^{-\top}\bar{\nabla}h)$, {i.e., geometric nonlinearities due to $\bar{F}$ are treated explicitly}. The implicit treatment of $e_z$ will provide the usual stabilisation of higher-order parabolic equations to avoid small time-step restrictions.
  
  \paragraph*{Step 2} Also for the discretized system, in the second step the ALE time-derivatives 
  $\dot{\bar{q}}=(\dot{\bar{\psi}},\dot{\bar{h}})$ are reconstructed from $\dot{h}$ computed in the first step. Therefore seek $(\dot{\bar{\psi}}:\bar\omega\to\R^d,\bar\lambda:\partial\bar\omega\to\R^d)$ such that
  \begin{align}
  &\int_{\bar\omega} \mathbb{D}\dot{\bar{\psi}}:\mathbb{D}\bar{v}_\psi\,\bar{J} {\rm d}\bar\omega + \int_{\partial\bar\omega}\bar{\lambda}\cdot\bar{v}_\psi\,{\rm d}\bar\gamma= 0,\\
  &\int_{\partial\bar\omega} \dot{\bar{\psi}}\cdot\bar{v}_\lambda \,{\rm d}\bar\gamma = \int_{\partial\bar\omega} -\bar{v}_\lambda\cdot\left(\frac{\bar{F}^{-\top}\bar{\nabla}\bar{h}}{|\bar{F}^{-\top}\bar{\nabla}\bar{h}|^2}\dot{h}+c_1 t\right)\,{\rm d\bar\gamma},
  \end{align}
  for all $(\bar{v}_\psi,\bar{v}_\lambda)$ with the Eulerian symmetric gradient $\mathbb{D}\bar{w}=\tfrac12(\bar{\nabla}\bar{w}\bar{F}^{-1}+\bar{F}^{-\top}\bar{\nabla}\bar{w}^\top)$. \\
  
  Now we discuss the choice of $c_1$ in more detail. If one is interested in gravity-driven droplet motion with gravity $\boldsymbol{g}=(g_x,g_z)$ and in-plane gravity $g_x=g\hat{e}$ pointing in the  direction of the unit vector $\hat{e}$, then often solutions contain a major translation coupled to a mostly unpredictable deformation of the overall droplet shape. For example: Neglecting the domain deformation and assuming only translation, the solution is a traveling wave of the form $h(t,\bx)=\bar{h}(\bar{\bx})$ for $\bar{\bx}=\bx-\hat{e}ct$ with velocity $c>0$. Such a profile has the Eulerian time-derivative
  $\dot{h}(t,\bx)=-c\hat{e}\cdot\bar\nabla\bar{h}(\bar{\bx})$ for all $\bx-\hat{e}ct=\bar{\bx}\in\bar\omega$. 
  If we would let $c_1\equiv 0$, then $\dot{\bar{\psi}}$ would only extract the normal component for the translation of $\omega$ and therefore according to \eqref{eqn:normalcomponents} we have $\dot{\bar{\psi}}=(\dot{h}/|\nabla h|)\nu=c(\hat{e}\cdot\nu)\nu$
  with outer normal $\nu$ of the translating support ($\bar{\nu}\equiv\nu$). In a co-moving frame we have the projected velocity 
  $\hat{e}\cdot(\dot{\bar{\psi}}-\hat{e}c)=c\bigl((\hat{e}\cdot{\nu}){\nu}-\hat{e}\bigr)\cdot\hat{e}=c(|\hat{e}\cdot{\nu}|^2-1)\le 0$. This shows that with normal updates, points will generally move backward in the comoving frame of reference, unless $\hat{e}\cdot{\nu}=\pm 1$.
  A good strategy in this case is to choose $c_1$ so that the domain motion can accomodate to the constant traveling wave speed. For an arbitrary given $\dot{\bar{\psi}}$, an estimate for the translational velocity $w=\hat{e}c$ can be made by minimizing
\begin{align}
\min_{w\in\R^d}\int_{\partial\omega} \left((\dot{\bar{\psi}}-w)\cdot\nu\right)^2\,{\rm d}\gamma \quad\rightarrow\quad w=\left(\int_{\partial\omega}\nu\otimes\nu{\rm d}\gamma\right)^{-1}\int_{\partial\omega}(\dot{\bar{\psi}}\cdot\nu)\nu{\rm d}\gamma\in\R^d,
\end{align}
and $\dot{\bar{\psi}}=(\dot{h}/|\nabla h|)\nu$. With the resulting estimated translational velocity we then choose $c_1=t(w\cdot t)$. Alternatively, one can also choose $\bar{v}_\lambda$ to point only in the direction of the Eulerian outer normal and determine the tangential velocity with the least elastic deformation energy.

\paragraph*{Step 3} The projection of the final step is performed with $\dot{h}$ and $\dot{\bar{\psi}}$ as in Definition~\ref{eqn:weakALEformulation} but instead the condition is required weakly, i.e., seek $\dot{\bar{h}}$ such that
\begin{align}
  \int_\omega \dot{\bar{h}}v\,{\rm d}\omega = \int_\omega (\dot{h}+\dot{\bar{\psi}}\cdot\nabla h)v\,{\rm d}\omega, 
\end{align}
in order to project from the possibly different finite element spaces for $\dot{h}\in V$ and $\dot{\bar{\psi}}\in W$ to the one for $\dot{\bar{h}}\in \bar{V}_0$.

{
\begin{remark}[Energy-stability]
   For gradient structures with 
  $\mathcal{R}(q,v):=\tau \Psi(q,v)+\mathcal{E}(q+\tau v)$,
  the successive minimization of $\mathcal{R}$ defines a sequence 
  \begin{align}
  \label{eqn:miniproblem}
  q^{n+1}=\mathop{\rm arg \min}_q \mathcal{R}(q^n,v), \qquad v=\frac1\tau(q-q^n).
  \end{align}
  Using the beforementioned properties $\Psi(q,0)=0$ and $\Psi(q,v)\ge 0$ the minimizer satisfies $\tau\Psi(q^n,v)+\mathcal{E}(q^{n+1})=\mathcal{R}(q^n,v)\le \mathcal{R}(q^n,0)\equiv \mathcal{E}(q^n)$ and therefore ensures the energy descent on the time-discrete level. Employing the same calculation as in the derivation of the saddle-point problem \eqref{eqn:derivation-saddle-point}, then the minimization of \eqref{eqn:miniproblem} gives
  \begin{align*}
    a(u,v) + \langle L(\tfrac{q^{n+1}-q^n}{\tau}),v\rangle &=0,\\
    \langle L(v),u\rangle &=\langle \mathrm{D}\mathcal{E}(q^{n+1}),v \rangle,
  \end{align*}
  where the state-dependence of $a(u,v)$ as in \eqref{eqn:miniproblem} is evaluated at $q^n$ and $\mathrm{D}\mathcal{E}$ is evaluated implicitly at $q^{n+1}$. For thin-film equations on a \emph{fixed domain} such a numerical scheme would be straightforward to implement. However, for the problem presented here the difference $\tfrac{q^{n+1}-q^n}{\tau}$ still needs to be transfered to a Lagrangian time-derivative of $h$ and the domain motion. Also, a fully implicit treatment of the dependence on the domain shape in $\mathcal{E}$ presents a highly nonlinear problem and we choose to use the decoupled semi-implicit scheme for its simplicity  and speed, since its formulated mostly in Eulerian variables and requires only the solution of a single linear system.
\end{remark}
}

\subsection{Time-stepping}
Based on the single decoupled step resulting in { unknown quantities that we denote by} $\dot{\bar{q}}_{\bar{q},\tau}=(\dot{\bar{\psi}},\dot{\bar{h}})$, different integration schemes to evolve $\bar{q}^{n-1}$ to $\bar{q}^n$ by a time-increment $\tau$ can be devised. Therefore, for a given $k\in\mathbb{N}$ we introduce sub-steps of size $\tau/k$
\begin{align}
  &\bar{q}_{(j)} = \bar{q}_{(j-1)}+\tfrac{\tau}k\dot{\bar{q}}_{\bar{q}_{(j-1)},\tau/k}, \qquad j=1\ldots k,
\end{align}
for 
$j=1\ldots k$, where we denote $\bar{q}_{(0)}=\bar{q}^{n-1}$ and $\bar{q}_{(k)}=\bar{q}^n_{\tau/k}$.
Then, the following higher-order strategies can be defined
\begin{align*}
  &\text{First-order semi-implicit base step (SEMI1):}\qquad&& \bar{q}^n(1,\tau) \equiv \bar{q}^n_\tau \equiv \bar{q}^{n-1}+\tau\dot{\bar{q}}_{\bar{q}^{n-1},\tau},\\
  &\text{Second-order Richardson extrapolation (RICH2):}\qquad&&
    \bar{q}^n(2,\tau) = 2 \bar{q}^n_{\tau/2} - \bar{q}^n_{\tau},\\
  &\text{Third-order Richardson extrapolation (RICH3):}\qquad&& 
    \bar{q}^n(3,\tau) = \tfrac{8}{3} \bar{q}^n_{\tau/4} - \tfrac{6}{3}\bar{q}^n_{\tau/2} + \tfrac{1}{3}\bar{q}^n_{\tau},
\end{align*} 
where for any $\bar{q}^n(r,\tau)$ of $r$th-order with step-size $\tau$, the higher-order method is defined by the extrapolation \cite{richardson1911approximate,constantinescu2010extrapolated}
\begin{align*}
  \bar{q}^n(r+1,\tau)=\frac{2^r\bar{q}^n(r,\tau/2)-\bar{q}^n(r,\tau)}{2^r-1},
\end{align*}
using the lower order methods and can thus be written in terms of sub-steps of the base step $\bar{q}^n_{\tau/k}$. The usefulness of these time-stepping schemes will be investigated with respect to the experimental convergence rate for exemplary solutions. Other higher-order discretisation schemes require different discretisations strategies of the original weak formulation, where instead of the semi-implicit decoupled step, usually the fully coupled formulation would be needed to discretized and solve. In practice this approach might be undesired because of the complexity of the resulting nonlinear PDE system. Note that while the overall problem is nonlinear, all elementary solution steps of the decoupled scheme involve solving linear systems once.

\subsection{Discretisation of strong contact line dissipation model}
Now we discuss the algorithmic aspects of the first-order in time (SEMI1) and higher-order space discretisation in the strong contact line dissipation limit and use the same higher-order time-discretisation based on extrapolation, i.e., RICH2, RICH3.
Therefore, for any given $q^{n-1}=(\omega^{n-1},h_{\omega^{n-1}})$ we determine the contact line velocity using \eqref{eqn:strongfriction_dynamic}, i.e., without the hats
\begin{align}
\dot{{x}}=-{n}[q^{n-1}]\big(-\tfrac\sigma2|\nabla h_{\omega^{n-1}}|^2+s\big).
\end{align}
Alternatively, as in \cite{grunewald2011variational}, this evolution can be regularized by adding line-tension $\epsilon\ge 0$ to the driving energy, i.e.,
\begin{align}
\label{eqn:stabilisation}
  \mathcal{E}(q)=\int_\omega e(x,h,\nabla h)\,{\rm d}\omega + \int_{\partial\omega} \epsilon \,{\rm d}\gamma,
\end{align}
which results in the regularized evolution of the boundary curve with the normal velocity
\begin{align}
\label{eqn:meancurvature}
  \dot{{x}}=-{n}(q^{n-1})\big(-\tfrac\sigma2|\nabla h_{\omega^{n-1}}|^2+s + \epsilon\kappa\big),
\end{align}
and $\kappa:\partial\omega\to\mathbb{R}$ being the mean curvature defined as $\kappa=-\nabla\cdot(\nabla h/|\nabla h|)$. For the stabilisation of quasi-static droplet evolution to be effective, an implicit time-discretisation of $\kappa$ using the Laplace-Beltrami technique as introduced by Dziuk \cite{dziuk1990algorithm} can be used. Note that the stabilisation in \eqref{eqn:stabilisation} can be applied to the hierarchy of models with contact line dissipation at any level. Using a vectorial $\mathbf{f}:\partial\omega\to\R^d$ with $\mathbf{f}=-n[q^{n-1}](-\tfrac\sigma2|\nabla h_{\omega^{n-1}}|^2+s)\nu$ we get the weak formulation for the velocity $\dot{x}=\dot{\mathbf{x}}\cdot\nu$
\begin{align}
\label{eqn:curveevolution}
\int_{\partial\omega} \dot{\mathbf{x}}\cdot \mathbf{v}+\tau\epsilon n \nabla_\Gamma\dot{\mathbf{x}}\cdot \nabla_\Gamma \mathbf{v} \,{\rm d}\gamma = \int_{\partial\omega} \mathbf{f}\cdot\mathbf{v}-\epsilon n \nabla_\Gamma\mathbf{x}\cdot \nabla_\Gamma \mathbf{v} \,{\rm d}\gamma,
\end{align}
that takes into account the mean curvature in \eqref{eqn:meancurvature} in a semi-implicit time discretisation.
This normal velocity needs to be subject to an extension to $\omega$ and a possible modification of its tangential component using the previously introduced ALE method. Then one can evolve the domain by setting $\omega^{n}=({\rm id}+\tau\dot{\Psi})(\omega^{n-1})$. Based on the updated domain, the new stationary droplet shape is determined by solving the minimization problem
\begin{subequations}
\label{eqn:statshape}
\begin{align}
    &\int_{\omega^n} {\pi} v_h\,{\rm d}\omega =\int_{\omega^n} e_z\cdot\nabla v_h + e_h v_h\,{\rm d}\omega, \\
   &\int_{\omega^n} {h}_{\omega^n}\,{\rm d}\omega =\mathit{vol}(q_{n-1}),
\end{align}
\end{subequations}
which gives the new state $q^n=(\omega^n,h_{\omega^n})$. While for the space-discretisation of $q^n$ we use higher-order isoparametric finite elements, the time-discretisation is semi-implicit. Since at every time  the stationary shape $h_{\omega}$ can be computed from the knowledge of $\mathit{vol}(q_{0})$ and the current domain shape $\omega$, this can be regarded as a nonlocal evolution for the boundary curve $\partial\omega$. 

\subsection{Implementation details}
The codes presented in the next section have been implemented both in \text{MATLAB} and \texttt{C++}. For the \texttt{C++} based implementation, we exploited the open-source library \texttt{deal.II}~\cite{ArndtBangerthBlais-2021-h, ArndtBangerthDavydov-2021-a}. The code supports arbitrary high order space discretisations, local mesh refinement, monolithic coupling of bulk and boundary terms, and geometrically preserving refinements and remeshing~\cite{HeltaiBangerthKronbichler-2021}.

We consider both scalar and vector $H^1$ conforming piecewise polynomial finite element spaces, based on a triangle discretisation (for the MATLAB implementation) and on a quadrilateral discretisation (for the \texttt{C++} implementation) $\bar\omega_h$ of the reference domain $\bar\omega$. In particular we consider the finite element spaces
\begin{equation}
    \text{P}k := \{ v \in C^0(\bar\omega_h) \text{ s.t. } v|_{T} \in \mathcal{P}^k(T) \qquad \forall T \in \bar\omega_h\},
\end{equation}
of piecewise polynomials of order $k$ on each triangle $T$ of $\bar\omega_h$ when considering triangular meshes, and
\begin{equation}
    \text{Q}k := \{ v \in C^0(\bar\omega_h) \text{ s.t. } v|_{T} \in \mathcal{Q}^k(Q) \qquad \forall Q \in \bar\omega_h\},
\end{equation}
of piecewise polynomials of order $k$ in each coordinate direction on each quadrilateral $Q$ of $\bar\omega_h$, when considering quadrilateral meshes.
The spaces defined on $\partial\bar\omega_h$ are constructed by restricting the spaces above. The deformed domain $\omega_h$ is constructed as the image of a field $\bar\psi \in \text{Q}k^2$ or in $\text{P}k^2$. We restrict all the examples to iso-parametric finite element discretisations, i.e., all fields are chosen in spaces with the same polynomial degree approximation $k$, ranging from one to three.

\section{Free boundary problems with dynamic contact angle: droplets and ridges}
\label{sec:results}
In this section we present results that underline the physical relevance of models with dynamic contact angle and show the convergence of the proposed spatial and temporal discretisations. We study the full transient model with Q1, Q2, and Q3 FE using decoupled time discretisation with SEMI1, RICH2 and RICH3 higher-order schemes. As an alternative for models with strong contact line dissipation, we also show related results with P1 and P2 FE using the same time discretisation.

Consider solutions $q(t)=\bigl(\omega(t),h(t)\bigr)$ of the thin-film problem with the energy 
\begin{align*}
    \mathcal{E}(q)=\int_\omega \tfrac\sigma2|\nabla h|^2 + s + \varrho h(g_x\cdot\mathbf{x}+\tfrac12 g_z h)\,{\rm d}\omega + \epsilon \int_{\partial\omega}{\rm d}\gamma,
\end{align*}
where we consider initial domains being balls of radius $R$ and arbitrary volume, i.e., $\omega(t=0)=B_R(0)$ and $\mathit{vol}(q_0)=V_0$. For the contact line dissipation we use a constant mobility $n[q]=n_0\in\mathbb{R}>0$. In order to nondimensionalize this problem we introduce scales
\begin{align*}
x=[L]\hat{x}, \quad z=[Z]\hat{z}, \quad t=[T]\hat{t},
\end{align*}
and rescale the initial data such that with $L=R$ we get $\hat{\omega}(t=0)=B_1(0)$ and with $Z$ such that $\mathit{vol}(\hat{q}_0)=1$. Using $\hat{e}_0=\sigma Z^2/L^2$ we can write
\begin{align*}
    \hat{\mathcal{E}}=(\hat{e}_0 L^d)^{-1}\mathcal{E}=\int_{\hat{\omega}} \tfrac12|\nabla \hat{h}|^2 + \hat{s} +  \hat{h}(\hat{g}_x\cdot\hat{\mathbf{x}}+\tfrac12 \hat{g}_z \hat{h})\,{\rm d}\hat{\omega} + \hat{\epsilon} \int_{\partial\hat{\omega}}{\rm d}\hat{\gamma},
\end{align*}
with $\hat{\epsilon}=\epsilon/(\hat{e}_0L)$, $\hat{s}=sL^2/(\sigma Z^2)$, $\hat{g}_x=\varrho g_x L^3/(\sigma Z)$ and $\hat{g}_z=\varrho g_z L^2/\sigma$, where for the validity of the thin-film approximation one would generally need to assume $Z\ll H$. While with these assumptions the initial data $\hat{\omega}(t=0)=B_1(0)$ and $\mathit{vol}(\hat{q}_0)=1$ are basically fixed, we have the three nondimensional parameters $\hat{s},\hat{g}_x,\hat{g}_z$ in the problem. Rescaling the time using $T=\mu L^4/(\sigma H^3)$ one obtains
\begin{subequations}
    \label{eqn:strong_thinfilm_rescaled}
    \begin{align}
      \partial_t \hat{h}-\nabla\cdot (m(\hat{h})\nabla\hat{\pi})=0,\qquad&\text{in }\hat{\omega},\\
      \partial_t \hat{h}+\hat{n}_0|\nabla \hat{h}|^2\hat{\zeta}=0,\qquad&\text{on }\partial\hat{\omega},\\
      m(\hat{h})\nu\cdot\nabla\hat{\pi}=0,\qquad&\text{on }\partial\hat{\omega},
    \end{align}
    together with the generalized forces $\hat{\pi}$ and $\hat{\zeta}$ being
    \begin{align}
    \hat{\pi}=-\Delta \hat{h} + \hat{g}_x\cdot\hat{\mathbf{x}}+ \hat{g}_z \hat{h},\qquad&\text{in }\hat{\omega},\\
    \hat{\zeta}=|\nabla \hat{h}|^{-1}\big(-\tfrac{1}{2}|\nabla \hat{h}|^2+\hat{s}+\hat{\epsilon}\hat{\kappa}\big),\qquad&\text{on }\partial\hat{\omega}.
    \end{align}
\end{subequations}
In the following we drop the hat symbols for convenience and use only the rescaled equations if not stated otherwise.

\subsection{Droplet motion with transient dynamics}
In this section we study solutions of the general transient model \eqref{eqn:strong_thinfilm_rescaled} with dynamic contact angle. The rescaled mobility in such a case is $m({h})=\tfrac13h^3+ bh^2$, where $b$ encodes the rescaled slip length. However, in order to fully concentrate on the effect of contact line dissipation and to avoid the discussion on the role of the slip length we will chose $m({h})={h}^2$ for this section. We consider the transient model without mean curvature stabilisation, i.e., $\epsilon\equiv 0$. 
\begin{figure}[H]
    \centering
    \includegraphics[width=0.95\textwidth]{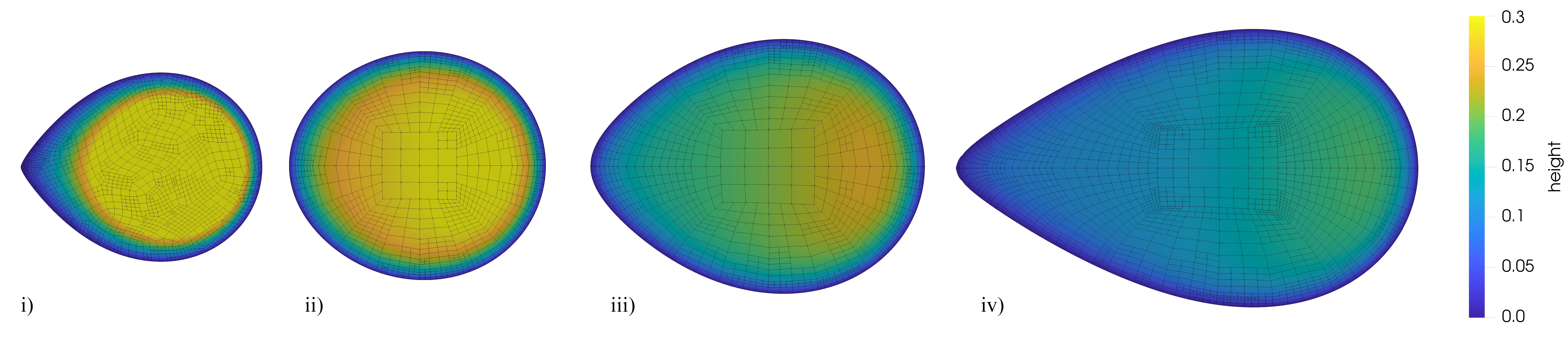}
    \caption{Sliding droplet traveling wave solution with Q2 isoparametric elements for $n_0=1$ and i) $g_x=5$, $g_z=0$, $s=1$ ii) $g_x=5$, $g_z=50$, $s=2$ iii) $g_x=5$, $g_z=50$, $s=1$ iv) $g_x=5$, $g_z=50$, $s=1/2$.} 
    \label{fig:slider_par}
\end{figure}  

In Fig.~\ref{fig:slider_par} traveling wave solutions for different parameters are shown. For sufficiently strong tangential gravity $g_x$ and small contact line mobility $n_0$ droplets develop a rather sharp corner as shown in i), whereas ii) and iii) show that increasing the normal gravity $g_z$ or the spreading coefficient $s$ can regularize this corner, whereas a smaller spreading coefficient $s$ will lead to more asymmetric droplets. The most notable feature that would be absent for $m(h)=h^2$ without dynamic contact angle is the sharp corner at the rear end of the sliding droplet \cite{peschka2015thin} but can otherwise only be recovered using cubic mobility $m(h)=\tfrac13 h^3+bh^2$ for sufficiently small regularisation $b\to 0$ \cite{schwartz2005shapes}. We are going to use the example of sliding droplets to study the convergence of solutions in space and time.

\begin{figure}[H]
    \centering
    \includegraphics[width=\textwidth]{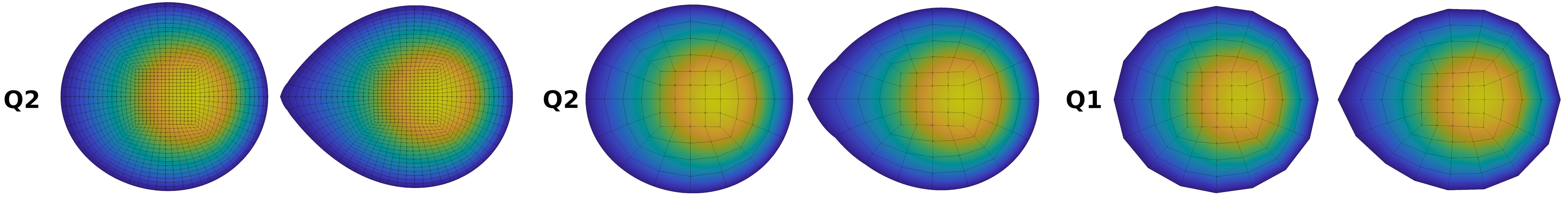}
    \caption{Sliding droplet (left) on fine mesh with isoparametric Q2 elements (middle) on a coarse mesh with isoparametric Q2 elements and (right) on a coarse mesh with isoparametric Q1 elements shown at $t=0.5$ and $t=2$ with ${g}_x=5$, ${g}_z=0$, ${s}=1$ and $\mathit{vol}({q}_0)=1$, $\tau=0.005$ and RICH2 time discretisation. The contact line dissipation is $n_0\equiv 1$.}  
    \label{fig:slider_space}
\end{figure}

In Fig.~\ref{fig:slider_space} traveling wave solutions with Q1 and Q2 FE on coarse and fine meshes are shown. The advantage of Q2 over Q1 elements is that droplets with a smooth circumference can be represented well already on coarse meshes as shown by the initial data at $t=0.5$. However, for pointy droplets at $t=2$ the higher order is not necessarily an advantage and local mesh refinement is the better strategy.

\begin{figure}[H] 
    \centering
    \includegraphics[width=0.4\textwidth]{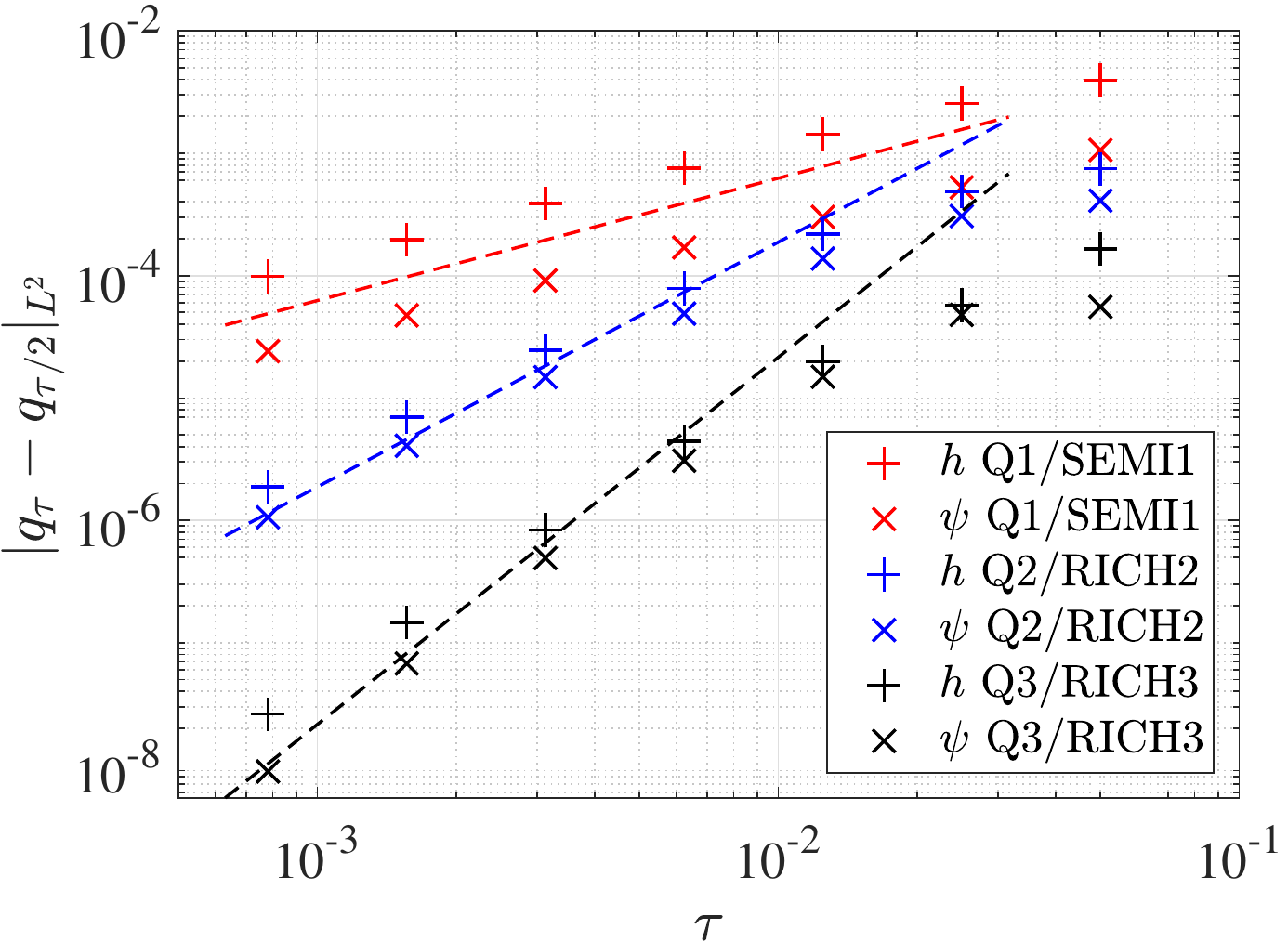}
    \caption{Convergence of $q=(\psi,h)$ in time with $g_x=2$, $g_z=0$, $n_0=\sigma=1$ at $T=0.1$ for Q1/SEMI1, Q2/RICH2, Q3/RICH3 space/time discretisation upon step-size bisection compared to linear (red dashed), quadratic (blue dashed), and cubic (black dashed) order.} 
    \label{fig:convergence_time} 
\end{figure}

The convergence of the higher-order time-discretisation is shown in Fig.~\ref{fig:convergence_time} and nicely produces the expected experimental order of convergence forthe $L^2$ errors of $h$ and $\psi$. However, the spatial order of convergence shown in Fig.~\ref{fig:convergence_time} only shows a marginal improvement using a higher-order discretisation. The main benefit for going from Q1 to Q2 appears to be a better prefactor for the error. This is not entirely unexpected since the regularity of weak solution, cf. \cite{giacomelli2013regularity,knupfer2015well}, does not give sufficient more regularity to justify a better convergence beyond the better representation of solutions of coarse meshes ({\color{blue}see \ref{appendix_A}}). Additionally, the projection of the Eulerian time-derivative $\dot{h}$ to the ALE time-derivative $\dot{\bar{h}}$ using $\dot{\bar{h}}=\dot{h}+\dot{\bar{\psi}}\cdot\bar{F}^{-\top}\bar{\nabla}\bar{h}$ introduces a nontrivial coupling of space and time discretisations, that could ultimately also affect the convergence order in space. Solving these problems requires special numerical techniques, which go somewhat beyond the scope of this paper, i.e., extension to geometrically conformal FE spaces or discontinuous Galerkin time discretisations cf.~\cite{boffi2004stability,bonito2013time}.

\begin{figure}[H]
    \centering
    \includegraphics[width=0.4\textwidth]{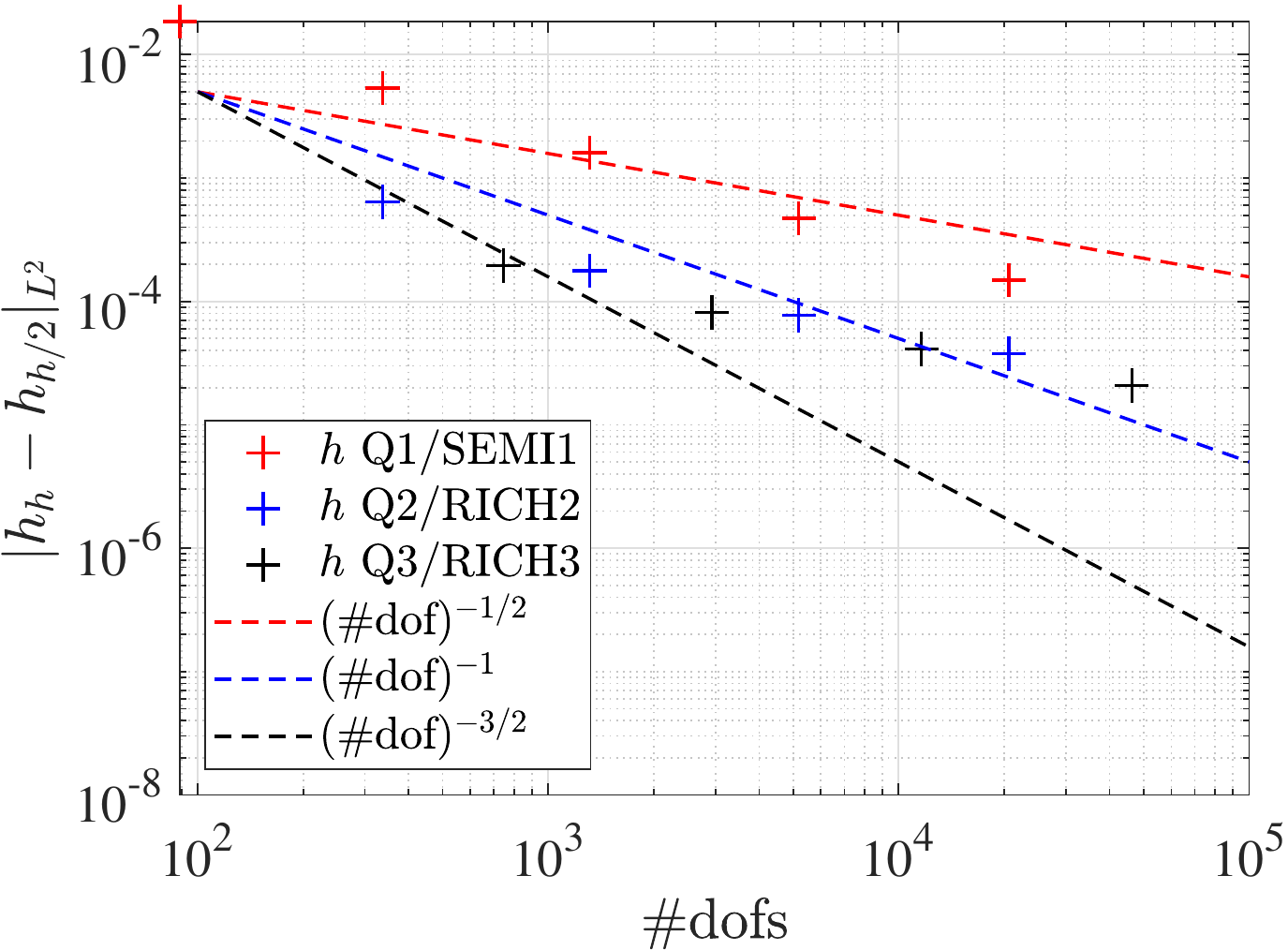}\hfill%
    \includegraphics[width=0.4\textwidth]{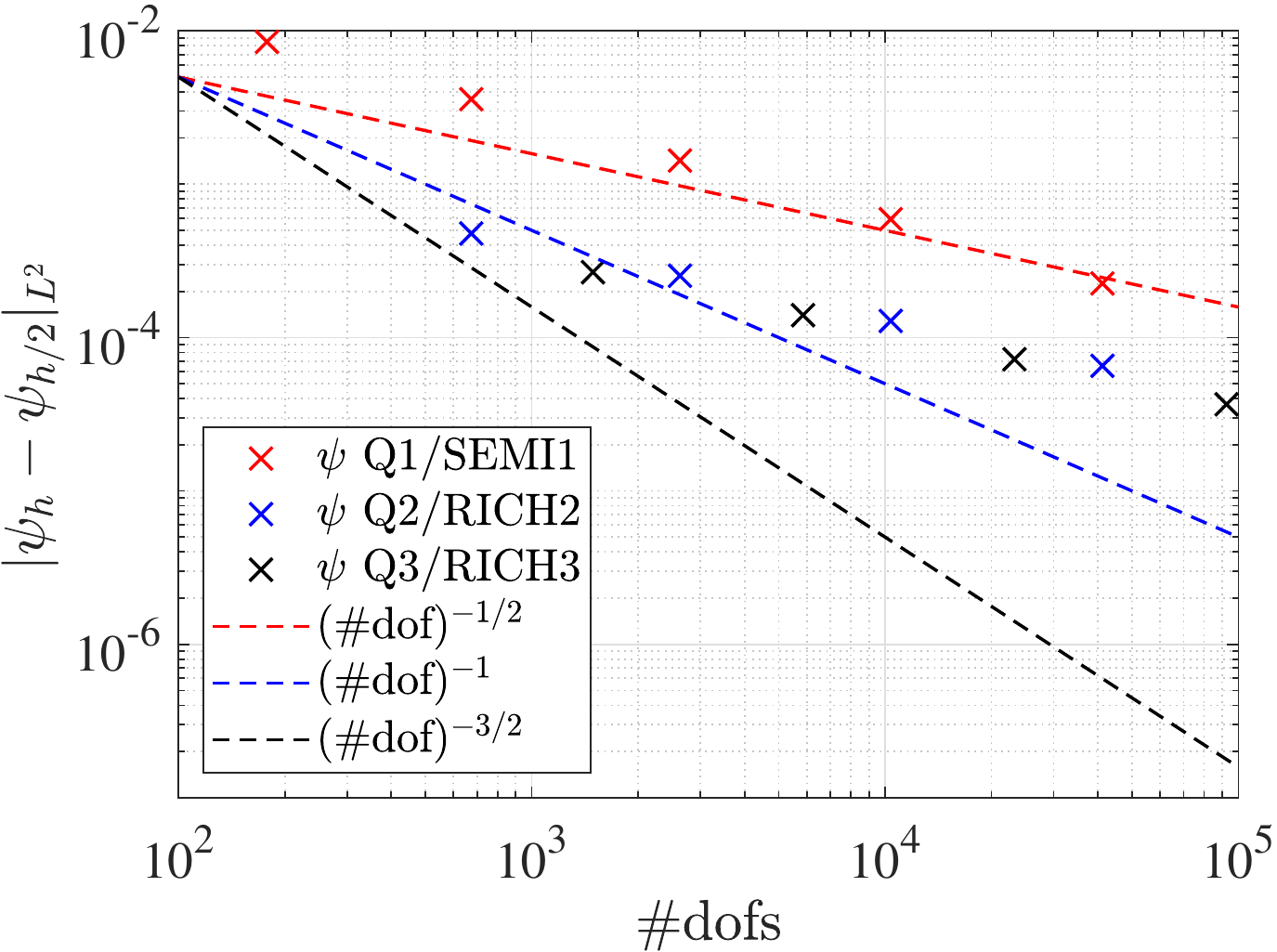}
    \caption{Convergence of $q=(\psi,h)$ in space with $g_x=2$, $g_z=0$, $n_0=\sigma=1,\tau=10^{-2}$ at $T=0.1$ for Q1/SEMI1, Q2/RICH2, Q3/RICH3 space/time discretisation upon uniform refinement for (left) $h$ and (right) $\psi$ using the $L^2$ norm compared to linear $h^1\sim(\mathrm{\#dof}^{-1/2})$ (red dashed), quadratic $h^2\sim(\mathrm{\#dof}^{-1})$ (blue dashed) and cubic $h^3\sim(\mathrm{\#dof}^{-3/2})$ (black dashed) order.} 
    \label{fig:convergence_space}  
\end{figure}  

\subsection{Droplet motion with strong contact line dissipation}

\begin{figure}[H]
\includegraphics[width=0.3\textwidth]{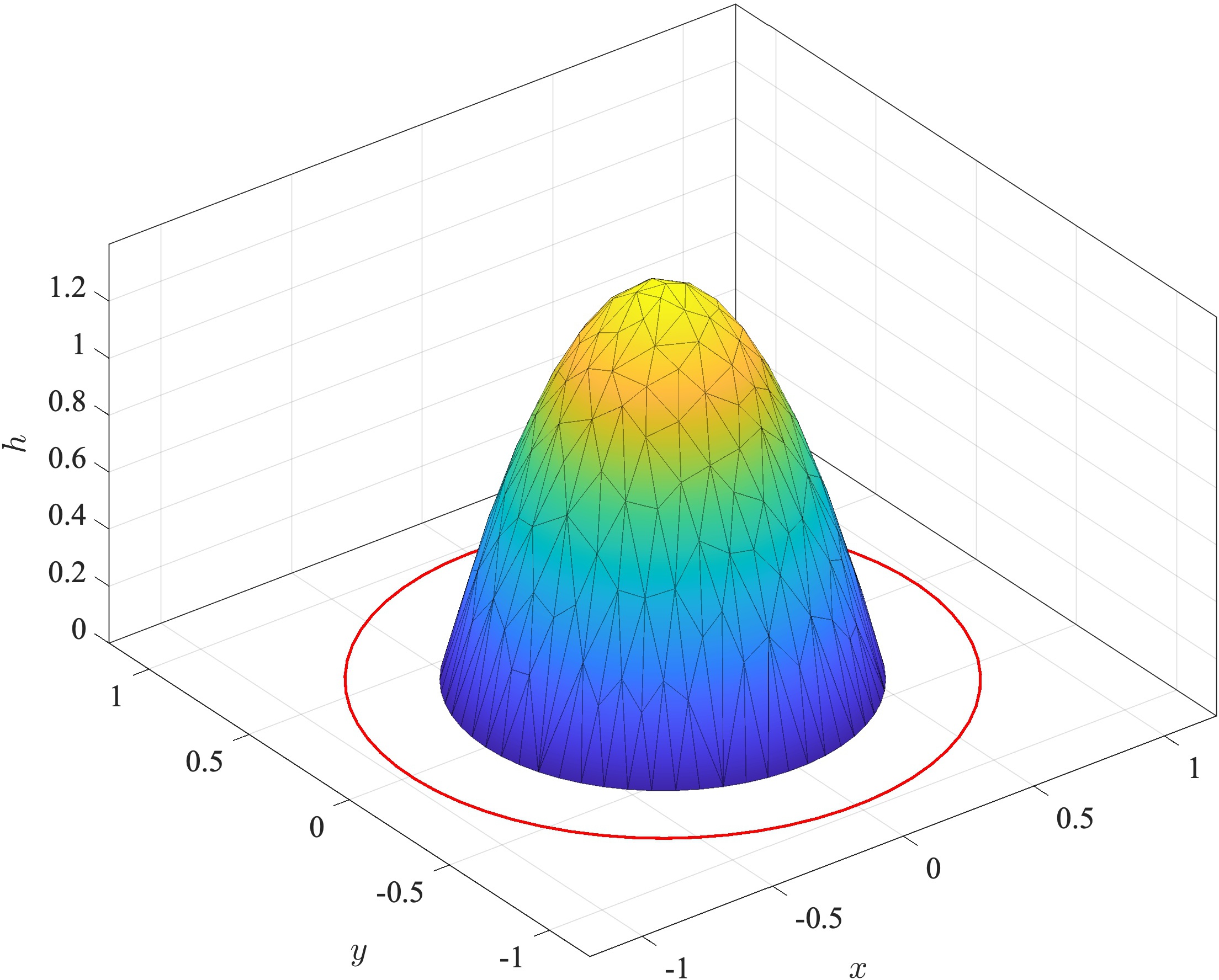}\hfill%
\includegraphics[width=0.3\textwidth]{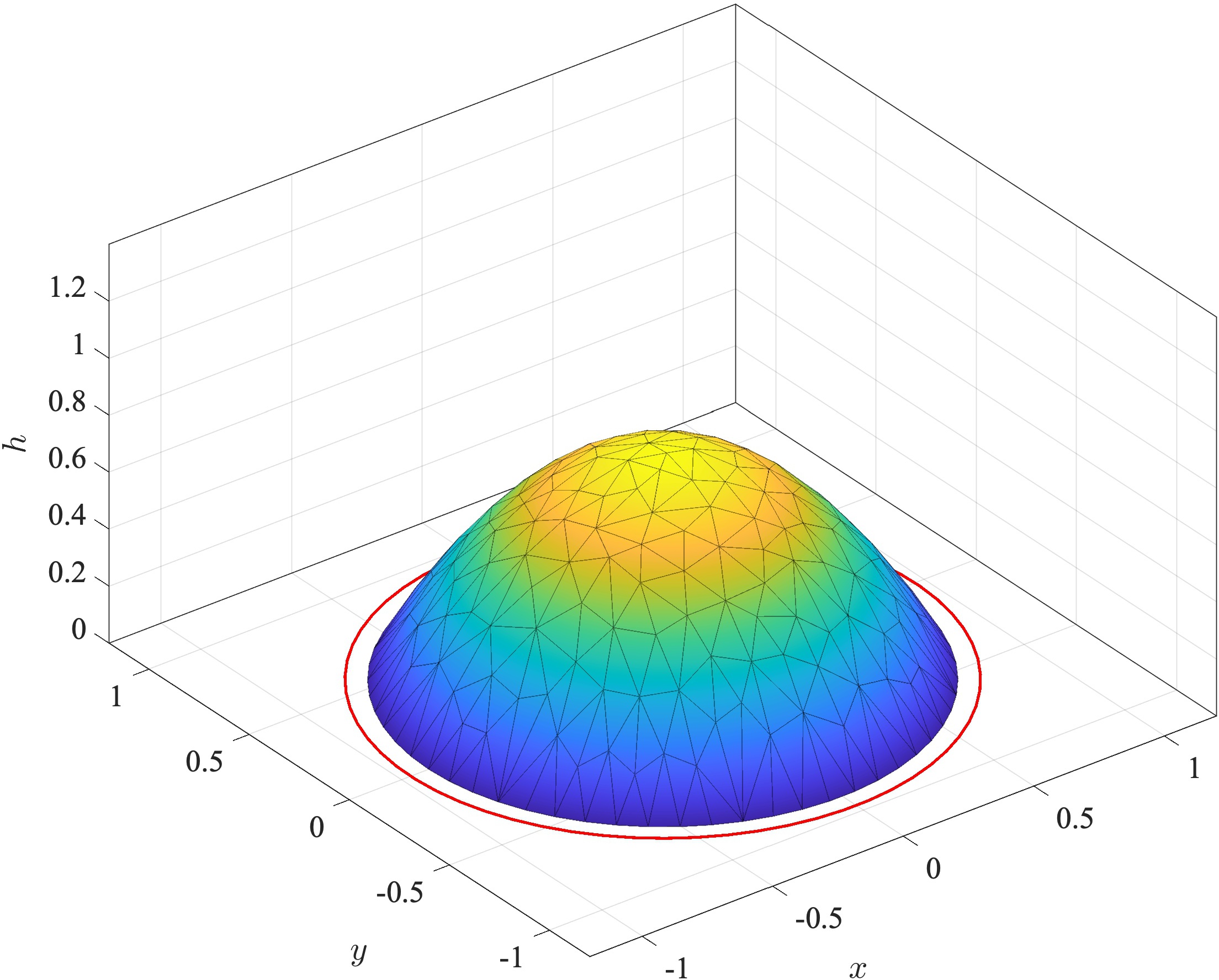}\hfill%
\includegraphics[width=0.3\textwidth]{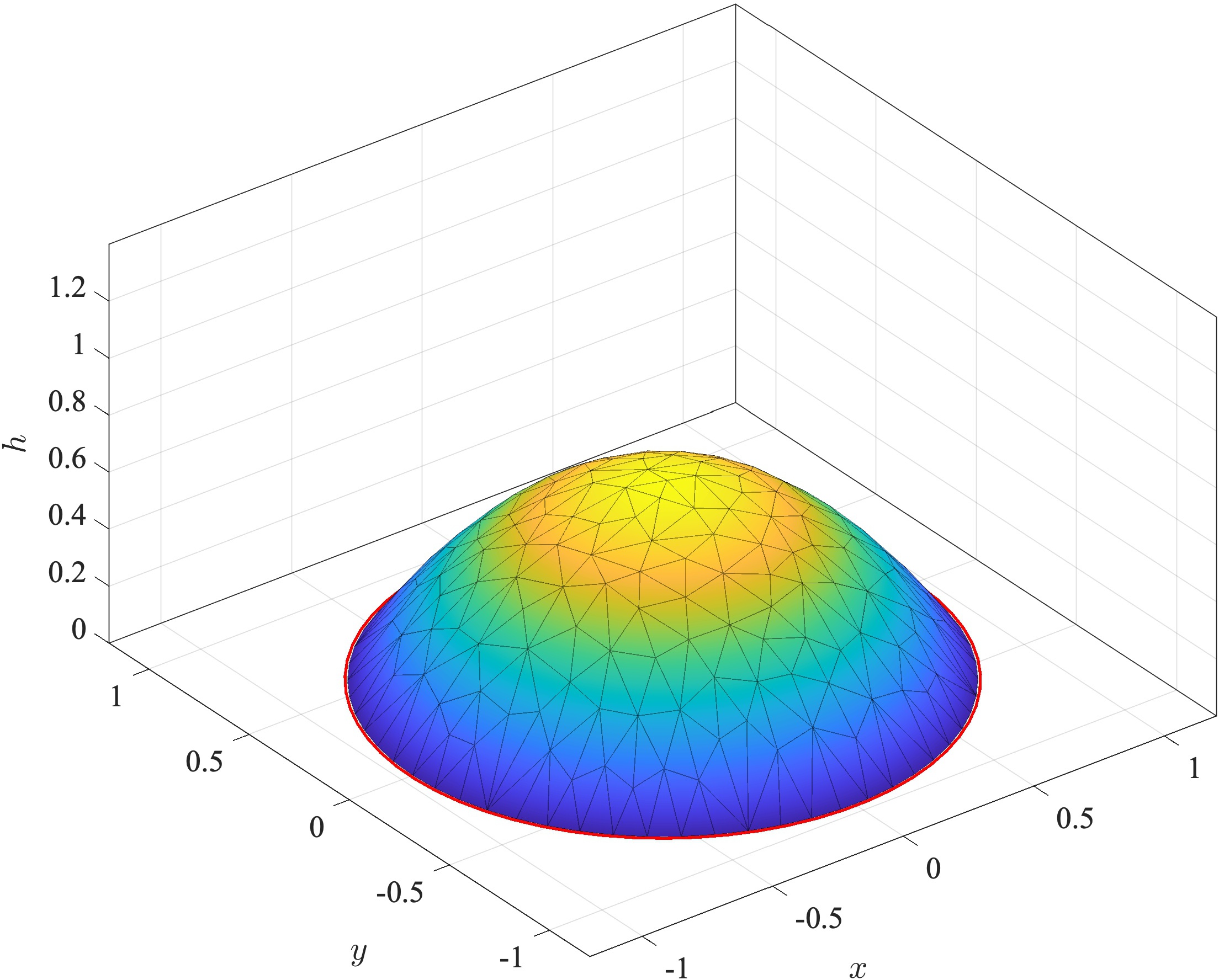}%
\caption{Stationary solutions of strong contact line friction model with ${s}=1$ and ${g}_x={g}_z=0$ obtained with $\tau=1/40$ and regularizing line tension ${\epsilon}=5.0,0.5,0.05$ decreasing from left to right in comparison with the radius of the droplet for ${\epsilon}=0$ shown as the red circle.}
\label{fig:statsol}
\end{figure}

When setting the time scale to $T=L/n_0{e}_0$, the nondimensional evolution for the dynamic contact angle becomes
\begin{align}
\label{eqn:normalvelocity}
    \dot{{x}}=-\left(-\frac12|\nabla{h}|^2+{s} + {\epsilon}{\kappa}\right).
\end{align} 
Without gravity ${g}_x,{g}_z=0$ the stationary state of this evolution for $d=2$ is given by ${\omega}_\textrm{stat}=B_r(0)$ and 
\begin{align}
    {h}_\textrm{stat}({\boldsymbol{x}})=c\left(1-\frac{|{\mathbf{x}}|^2}{r^2}\right), \qquad r=\left(\frac{4}{\pi\sqrt{2{s}}}\right)^{1/3},\qquad c=\frac{2}{\pi r^2},
\end{align}
provided that the regularisation is sufficiently weak, i.e., ${\epsilon}\ll {s}r$. The following examples are solved with P2 isoparametric finite elements and with the semi-implicit time-discretisation \eqref{eqn:statshape} and \eqref{eqn:curveevolution}. If one is interested only in the long-time behavior of solutions, then the extra degree of freedom in the scaling can be used to set ${s}=1$. Stationary solutions for different ${\epsilon}$ are shown in Fig.~\ref{fig:statsol} and agree with the exact solution for ${\epsilon}\to 0$.

\begin{figure}[H]
    \centering
    \includegraphics[width=0.9\textwidth]{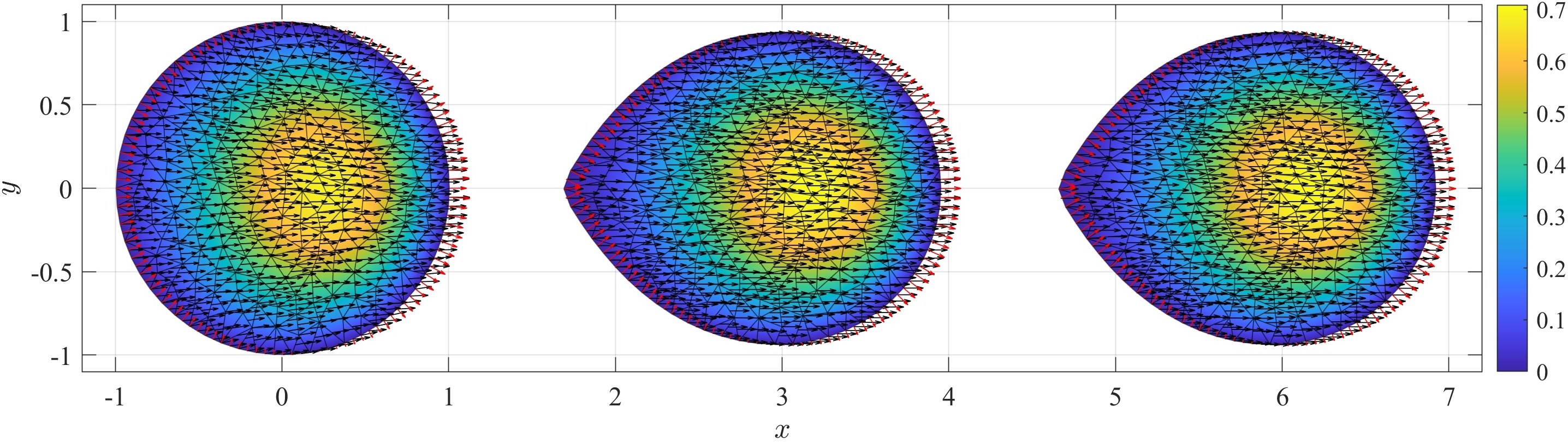}%
    \caption{Sliding droplet with ${s}=1$, ${g}_x=13/4$, ${g}_z=0$, ${\epsilon}=0.02$ at times $t=0,5/2,5$ from left to right. The shading encodes the height ${h}$, the normal velocity on the boundary shows \eqref{eqn:normalvelocity} (red arrows) and the reconstructed ALE velocity field $\dot{\bar{\psi}}$ is depicted with black arrows.}
    \label{fig:twsol}
\end{figure}

\begin{figure}[H]
    \centering
    \includegraphics[width=0.7\textwidth]{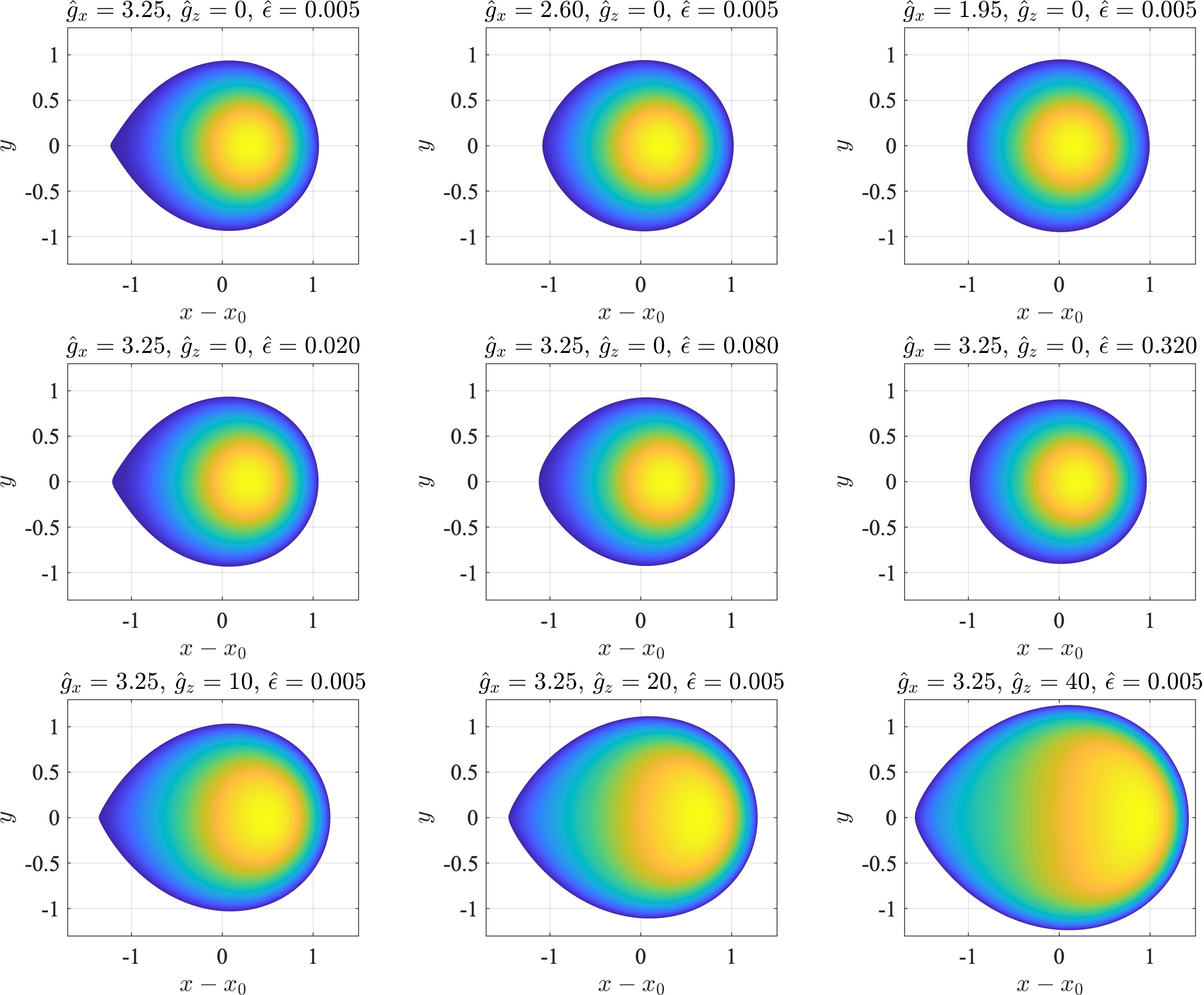}%
    \caption{Sliding droplet with different parameters ${g}_x,{g}_z,{\epsilon}$ and all with ${s}=1$. Positive solution exist for a certain time interval or traveling wave solutions exist for all times for sufficiently small ${g}_x < g_\textrm{crit}$. For larger ${g}_x>g_\textrm{crit}$ solutions become eventually become negative and loose their meaning. At a ${g}_x=g_\textrm{crit}$ the solution develops a trailing corner that can be regularized using the line tension ${\epsilon}$.}
    \label{fig:twsolparameters}
\end{figure}

With nonzero but sufficiently small tangential gravity $|{g}_x|$ we obtain traveling wave solutions, i.e., height profiles that move with constant translational velocity as in the example shown in Fig.~\ref{fig:twsol}. The corresponding vector fields shown for the solution in Fig.~\ref{fig:twsol} prove the effectiveness of the ALE method, i.e., an almost constant translational vector field is recovered from the normal boundary velocity field. In Fig.~\ref{fig:twsolparameters} we show the parameter space of different traveling wave solutions with dynamic contact angle with different horizontal gravity, normal gravity, and mean curvature regularisation. Again, one can observe that while stronger horizontal gravity leads to droplets with a corner at the trailing edge, normal gravity and mean curvature tend to reduce this curvature again.

\begin{figure}[ht]
    \centering
    \includegraphics[width=0.4\textwidth]{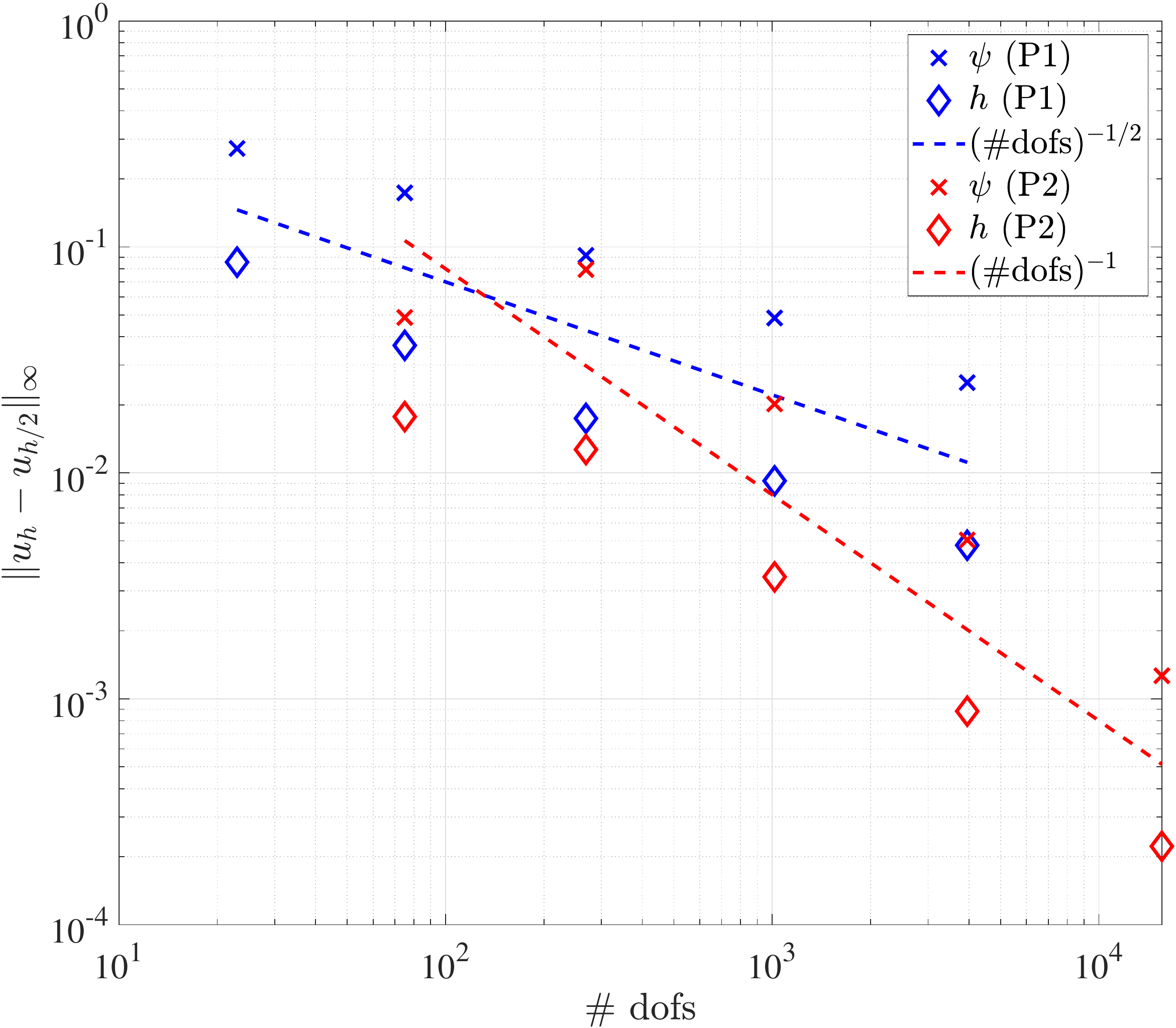}%
    \includegraphics[width=0.4\textwidth]{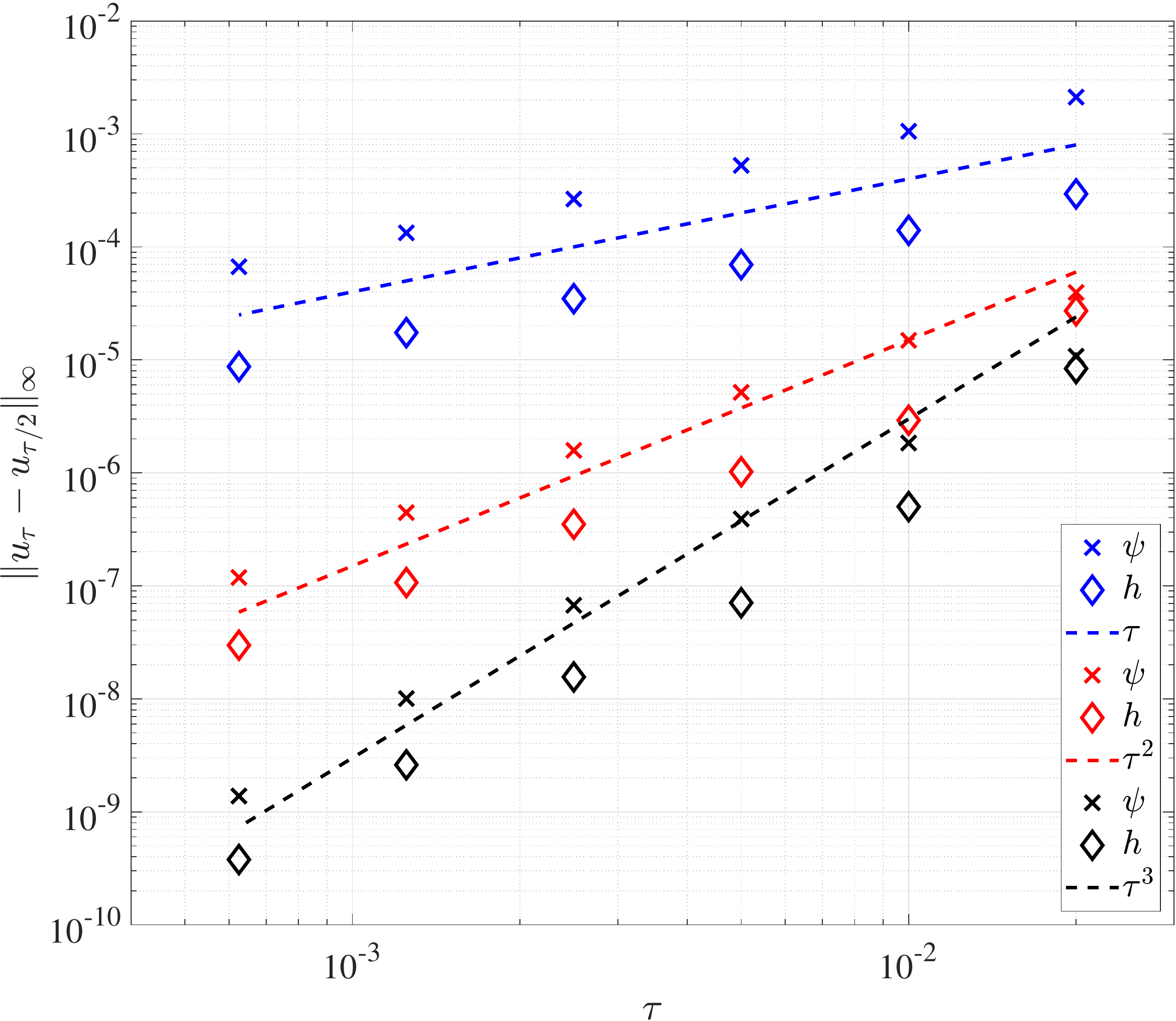}%
    \caption{Experimental order of convergence in max-norm with strong contact line dissipation (left) in space using uniform refinement for P1 and P2 finite elements with for $T=0.8$ and $\tau=0.002$ with RICH3 time discretisation and (right) in time using stepsize bisection for SEMI1, RICH2, RICH3 with $T=1$ and P2 finite elements and the number of degrees of freedom for $(\psi,h)$ is $\#\mathrm{dofs}=((863,863),863)$.}
    \label{fig:twerrors}
\end{figure}

Most importantly, the convergence in space and time as shown in Fig.~\ref{fig:twerrors} shows the expected higher-order convergence both in space and time. Again, we can achieve higher-order convergence in time using the Richardson extrapolation RICH2 (red) or RICH3 (black) giving second and third order convergence in $\tau$ compared to the linear order of the SEMI1 (blue) time-discretisation. But now we also achieve a higher-order convergence in space using P2 isoparametric elements, which is expected since the regularity of the energy-minimizers is better than the one of traveling wave solutions of the full transient problem.

\subsection{Instability of a ridge}

In this section we are finally going to consider the problem of a liquid ridge, which becomes unstable and pinches off~\cite{king2009linear}. Therefore consider the initial rectangular domain $\bar\omega=[0,L]\times[0,H]\subset\mathbb{R}^2$ where we have the free boundary $\bar\gamma_{\rm f}=\{0\}\times [0,H] \cup \{L\}\times[0,H]$ and the sliding boundary $\bar\gamma_{\rm s}=\partial\bar{\omega}\setminus\bar\gamma_{\rm f}$. We deform this domain and its boundaries to $\omega\equiv\omega(t=0)$ with $\bar\psi:\bar\omega\to\omega$ defined as 
\begin{align}\bar\psi(\bar{x},\bar{y})=\begin{pmatrix}\bar{x}\\\bar{y}\end{pmatrix}+ \delta\cos(2\pi\bar{y}/H)\begin{pmatrix}\bar{x}\\0\end{pmatrix},
\end{align}
and define the state space $Q=\{(\psi,h):h|_{\gamma_{\rm f}}=0,\psi|_{\gamma_{\rm s}}\cdot e_y=0\}$, i.e., $\psi:\omega_0\to\mathbb{R}^2$ only supports motion in $x$-direction on $\gamma_{\rm s}$ and $h:\omega_0\to\mathbb{R}$ vanishes on $\gamma_{\rm f}$. The main difference compared to the droplet dynamics is that we have introduced a sliding boundary, on which $h$ needs not to vanish but can slide. For these functions $q\in Q$ we minimize the energy $\mathcal{E}(q)$ subject to the constraint of a given volume $\mathit{vol}(q)$ and evolve $q(t)$ according to \eqref{eqn:meancurvature}. In comparison to the droplet motion we will consider two cases, i.e., the transient dynamics with $n[q]=n_0$ and the strong contact line friction limit with $n[q]=|\nabla h|^\theta$ for some $\theta\in\{-1,0,+1\}$.

\begin{figure}[H]
    \centering
    \includegraphics[height=0.22\textwidth]{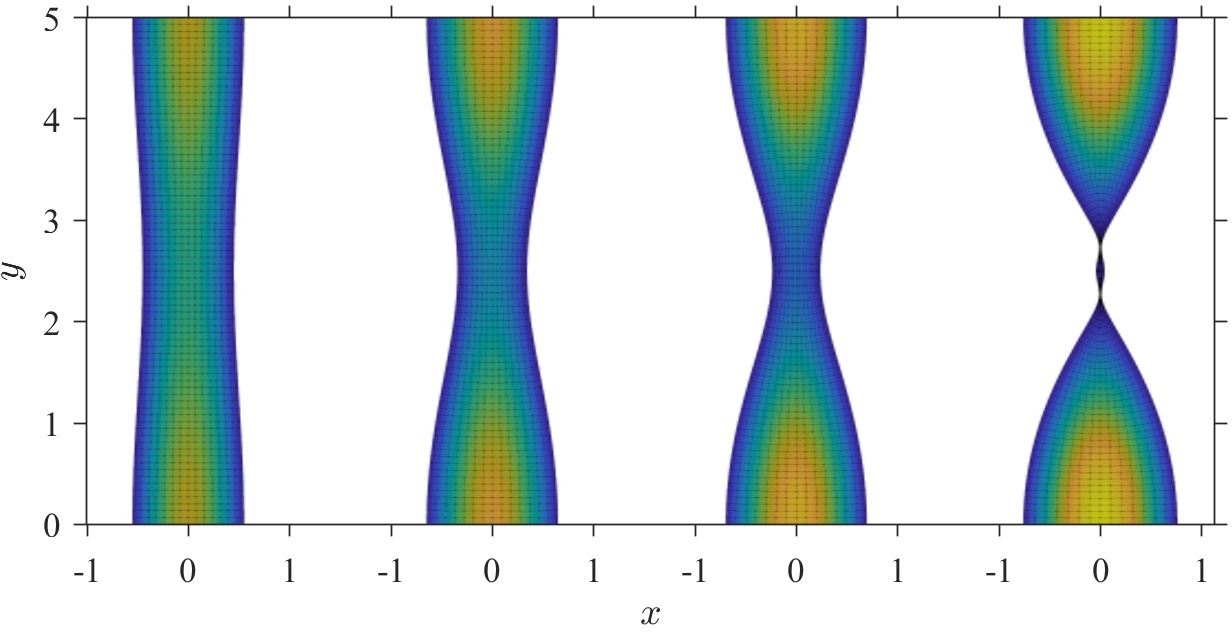} 
    \includegraphics[height=0.22\textwidth]{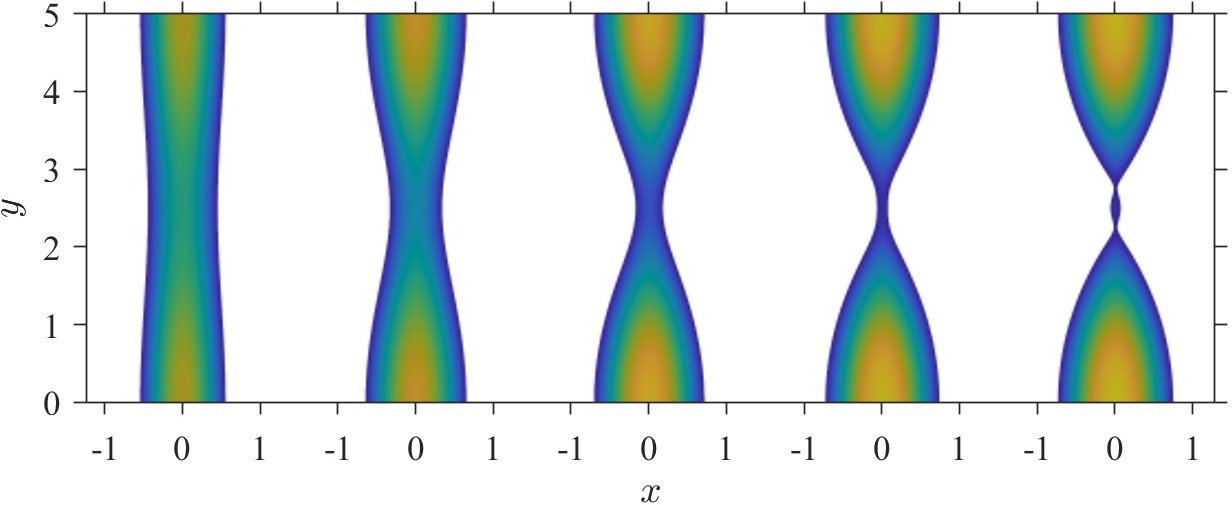}
    \includegraphics[height=0.22\textwidth]{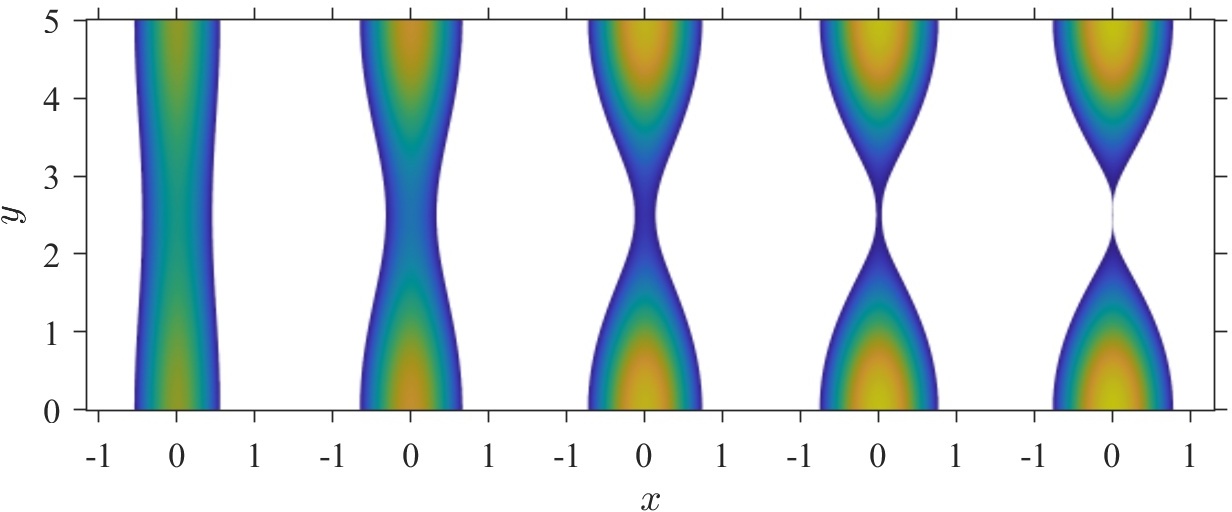} 
    \caption{Instability of a ridge with transient dynamic contact angle with $s=1$, $g_x=g_z=0$ and $\epsilon=0$ and different contact line friction, i.e.,  (top left) $n_0=1$ (top right) $n_0=10$ and (bottom) $n_0=0.1$ for mobility $m(h)=h^2$ at different times progressing from left to right.}
    \label{fig:ridge_dyn}
\end{figure}

The dynamics of the ridge using the transient dynamics with $n_0\in\{0.1,1,10\}$ is shown in Fig.~\ref{fig:ridge_dyn}. For both $n_0=1$ and $n_0=10$ the pinch-off does not occur symmetrically at $y=H/2$ but rather forms a satellite droplet. A similar effect was observed with thin-film type models with precursor in \cite{peschka2019signatures} and was triggered by different mobilities, i.e., $m(h)=h^2$ versus $m(h)=h^3$. This behavior changes for $n_0=10$, where the rupture point moves to the center and a satellite droplet is not visible anymore. Close inspection near the pinch-off still shows a much smaller remaining droplet.

\begin{figure}[ht]
    \centering
    \includegraphics[height=0.2\textwidth]{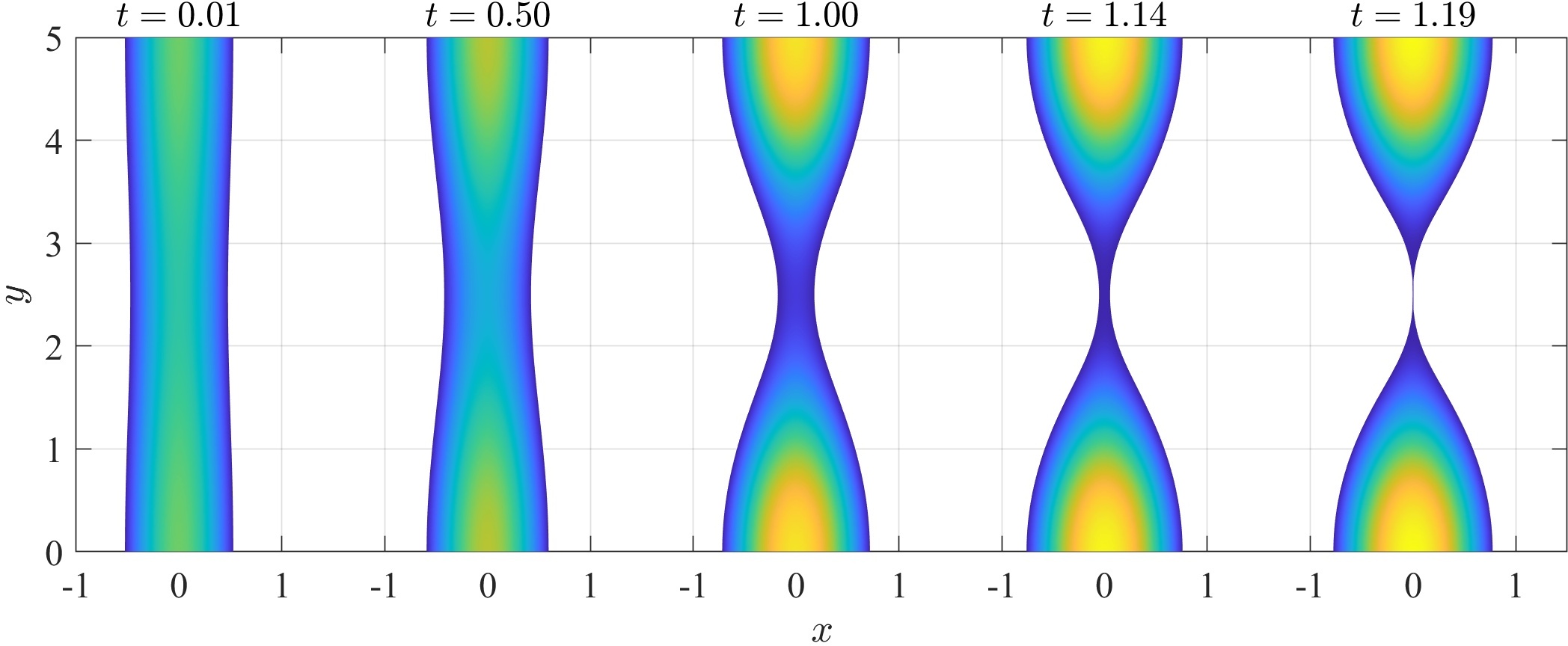} 
    \includegraphics[height=0.2\textwidth]{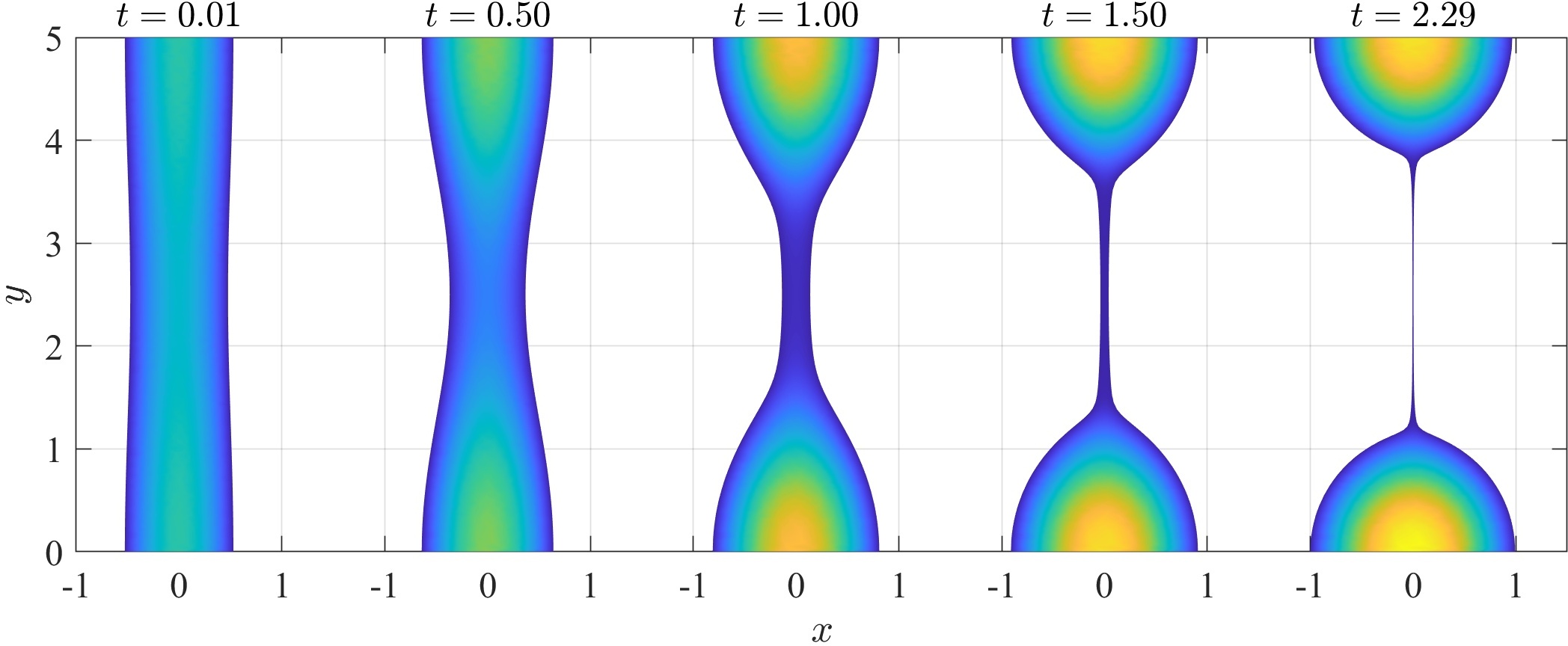}
    \includegraphics[height=0.2\textwidth]{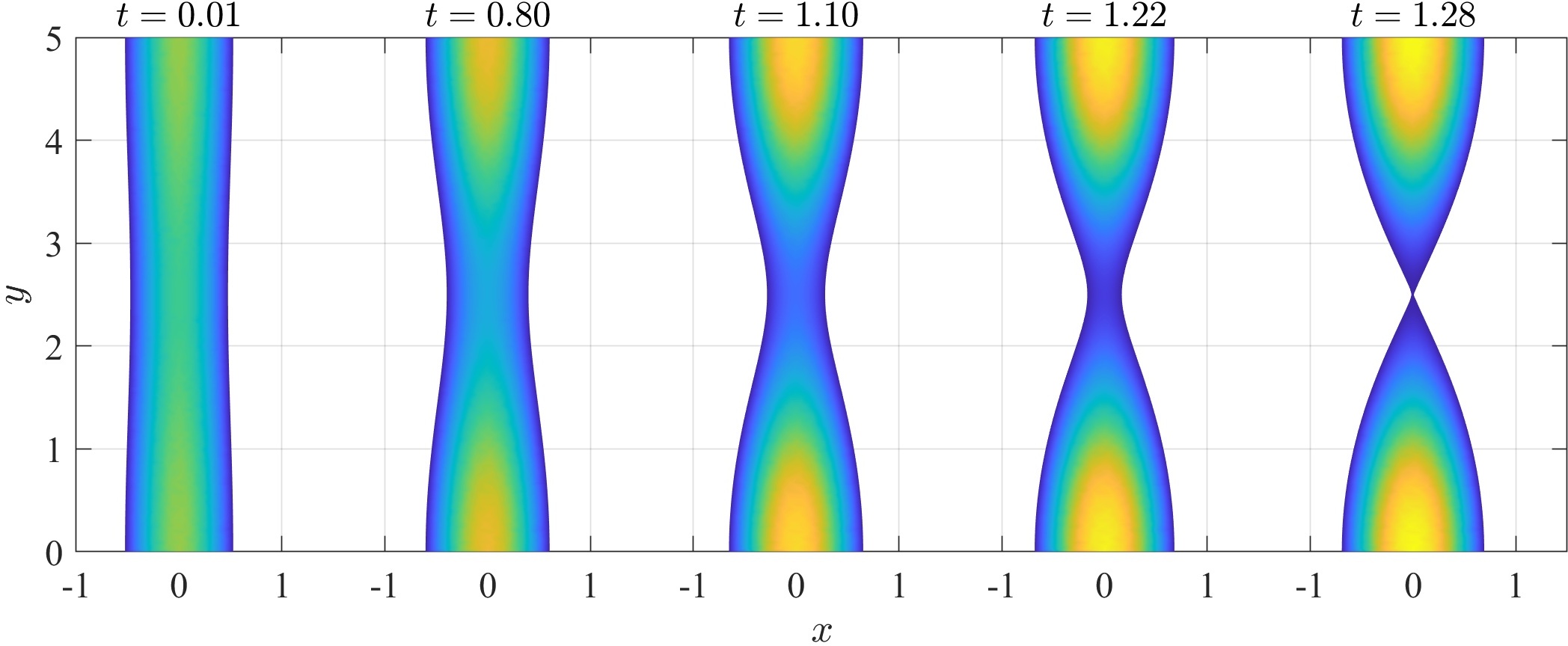} 
    \caption{Instability of a ridge in strong contact line dissipation model with $s=1$, $g_x=g_z=0$ and $\epsilon=0.02$ and different laws for the dynamic contact angle, i.e.,  (top left) $n[q]=n_0$ (top right) $n[q]=n_0|\nabla h|$ and (bottom) $n[q]=n_0/|\nabla h|$ at different times for $n_0=1$.}
    \label{fig:ridge_strong}
\end{figure}

The corresponding ridge instability in the strong contact line dissipation limit is shown in Fig.~\ref{fig:ridge_strong}. The solution shown in the upper row has a pinch-off at a finite time, i.e., the ridge width goes to zero as $w\sim(t_c-t)$ for some finite time $t_c$. This solution corresponds to the limit $n_0\to 0$ of the transient model and has the same qualitative ridge shape at the pinch-off as the transient solution with $n_0=0.1$ but without satellite droplets. Due to the scaling invariance it does not make sense to study different values $n_0$ in the strong dissipation limit. However, a simple asymptotic analysis shows that for $n=n_0$ the contact angle approaches zero at the pinch-off and thus choosing $n[q]=|\nabla h|^\theta$ with different $\theta$ should be expected to have a drastic effect on the instability. Indeed, for $\theta=1$ we observe the formation of two almost stationary droplets connected by a thin thread of liquid with thickness of the thread going to zero exponentially (in infinite time), i.e., $w\sim \exp(-ct)$ for some constant $c$.
For $\theta=-1$ we observe a finite time pinchoff with $w\sim (t_c-t)^{1/2}$ for the finite rupture time $t_c$. 

The examples presented here underline the relevance of correctly capturing the physics of dynamic contact angle for wetting flows, since the resulting shapes of drops and ridges may differ drastically. 

\section{Conclusion}
We have shown that the thin-film free boundary problem with moving contact lines posesses an energetic-variational structure, which allows to write the dynamics for $q=(\psi,h)$ in the form $q=-K(q){\rm D}\mathcal{E}(q)$, where the $K(q)$ operator contains an appropriate coupling of bulk and interface terms leading to dynamic boundary conditions for the dynamic contact angle. In the limit of strong dissipation at the contact line we derived a simpler operator $K$, where at each moment $h$ minimizes the energy with fixed $\psi$ using a given volume $\mathit{vol}(q)=\int h\,{\rm d}x$.

We proposed space and time discretisations for these gradient structures and also implement corresponding ALE techniques for the motion of the domain using $\psi$. For the transient problem we proposed a simple decoupled discretisation, which keeps the solution strategy for the free boundary problem simple, i.e., solution of a weak formulation, determination of boundary velocities, and ALE reconstruction and update of the solution.
For the strong contact line dissipation limit we added line tension as means of regularisation for $\psi$ and otherwise use a similar discretisation strategy.
In both cases we established higher-order discretisation strategies in space and time and performed convergence tests for the transient model and the strong dissipation model.

Thin-film type models with precursors rely on mobilities $m(h)=\tfrac13h^3+bh^2$ with sufficiently small $b$ to generate a logarithmic singularity and thereby localize dissipation near the contact line. In practice, the usage of a precursor presents another regularisation parameter and potentially renders smaller values for the contact line dissipation $n_0$ difficult to attain and thus, the practical impact of contact line dissipation of dewetting flows remains mainly unexplored beyond the 1D setting, where mostly exact solutions (self-similar or traveling wave solution) or asymptotic expansions have provided some understanding in the past. Our results show how the consistent use of dynamic contact angle in a higher-dimensional setting can potentially generate new types of solutions, e.g., pinch-off for different choices of $n[q]$.

\section*{Acknowledgement}

DP thanks Manuel Gnann, Lorenzo Giacomelli and Leonie Schmeller for fruitful discussions on the topic and acknowledges the financial support within the DFG-Priority Programme 2171 \emph{Dynamic Wetting of Flexible, Adaptive, and Switchable Substrates} by project \#422792530. LH acknowledges support from the National Research Project (PRIN  2017) ``Numerical Analysis for Full and Reduced Order Methods for the efficient and accurate solution of complex systems governed by Partial Differential Equations'', funded by the Italian Ministry of Education, University, and Research. Both authors thank the Berlin Mathematics Research Center MATH$^+$ and the  Einstein Foundation Berlin for the financial support within the Thematic Einstein Semester \emph{Energy-based mathematical methods for reactive multiphase flows} and by project AA2-9.

\newpage
\appendix
{
\section{Verification of expected order of convergence}
\label{appendix_A}

While the convergence order for the transient thin-film model in Fig.~\ref{fig:convergence_space} might surprise at first, this is related to the lack of regularity of solutions and therefore expected. For complete wetting, in \cite{giacomelli2013regularity} it is shown that source-type solutions are of the form $h(t,x)=t^{-\hat{\alpha}} H(xt^{-\hat{\alpha}})$ with $H(\xi)\sim \xi^{\hat{\nu}} (1+v(\xi,\xi^{\hat{\beta}}))$ for some rational exponents $\hat{\alpha},\hat{\nu}$ and irrational exponent $\hat{\beta}$. Similarly, in the partial wetting regime general solution of the thin-film problem feature a traveling-wave-type solution near a moving contact line, which in general prevents solutions to be smooth at the contact line. The exact singular behavior depends on the degeneracy of the elliptic operator $L(h)u = -\nabla\cdot (m(h)\nabla u)$ that appears in the thin-film problem. Here we have  $m(h)\sim h^2 \sim \mathrm{dist}(x,\partial\omega)^2\to 0$. 

For example, assume $\omega=(0,1)\subset\R$, then solutions $u(x)$ of $Lu(x)=f(x)$ with $L u(x) =-\nabla \cdot (\mu(x) \nabla u(x)) + u(x)$ for $\mu(x)=x^2$ have the same degeneracy as the thin-film problem on the left boundary. For $f(x)=x$ and with natural boundary condition we have the exact solution
\begin{align}
  u(x)=\frac{1}{c}x^c - x, \qquad c=\frac{\sqrt{5}-1}{2} = 0.6180...,
\end{align}
which lacks smoothness at $x=0$, i.e., $u''(x)$ is not square integrable. On the other hand, using $\mu(x)=1+x^2$ entirely resolves this issue and solutions converge with their optimal convergence order. In Fig.~\ref{fig:suboptimal_convergence} we see suboptimal $L^2$ convergence with $\mu=x^2$, which remains first order but improves the prefactor for higher P$n$. On the other hand, for $\mu(x)=1+x^2$ we see the optimal convergence order $\mathcal{O}(h^{n+1})$ for P$n$ elements and mesh size $h$. We claim that, in essence, this is the behavior seen for the transient solution of the thin-film problem in Fig.~\ref{fig:convergence_space}. On the other hand, the solution of the strong contact-line dissipation limit is smooth and therefore converges optimally, see Fig.~\ref{fig:twerrors}. 

\begin{figure}[H]
\centering
  \includegraphics[width=0.48\textwidth]{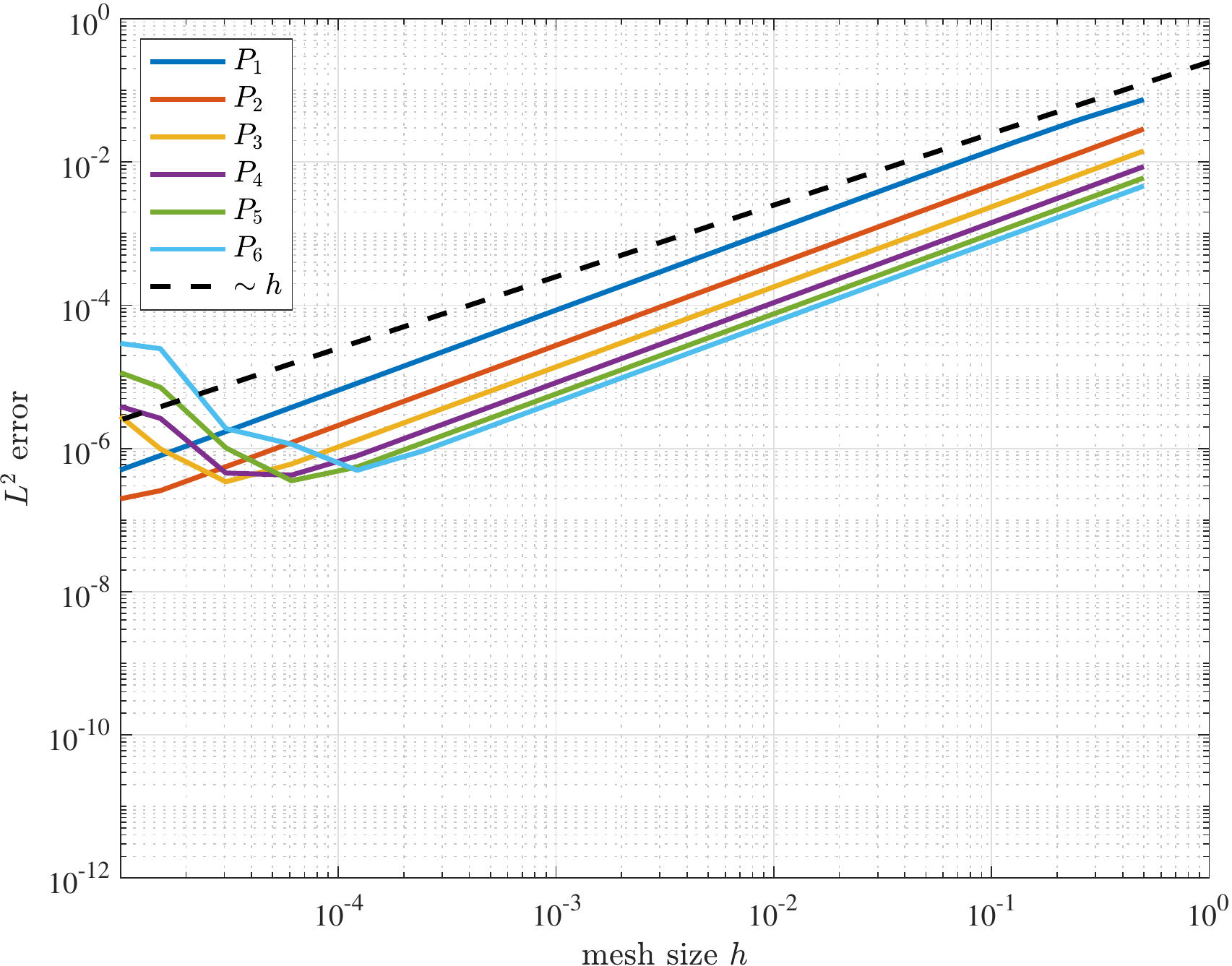}\hfill
  \includegraphics[width=0.48\textwidth]{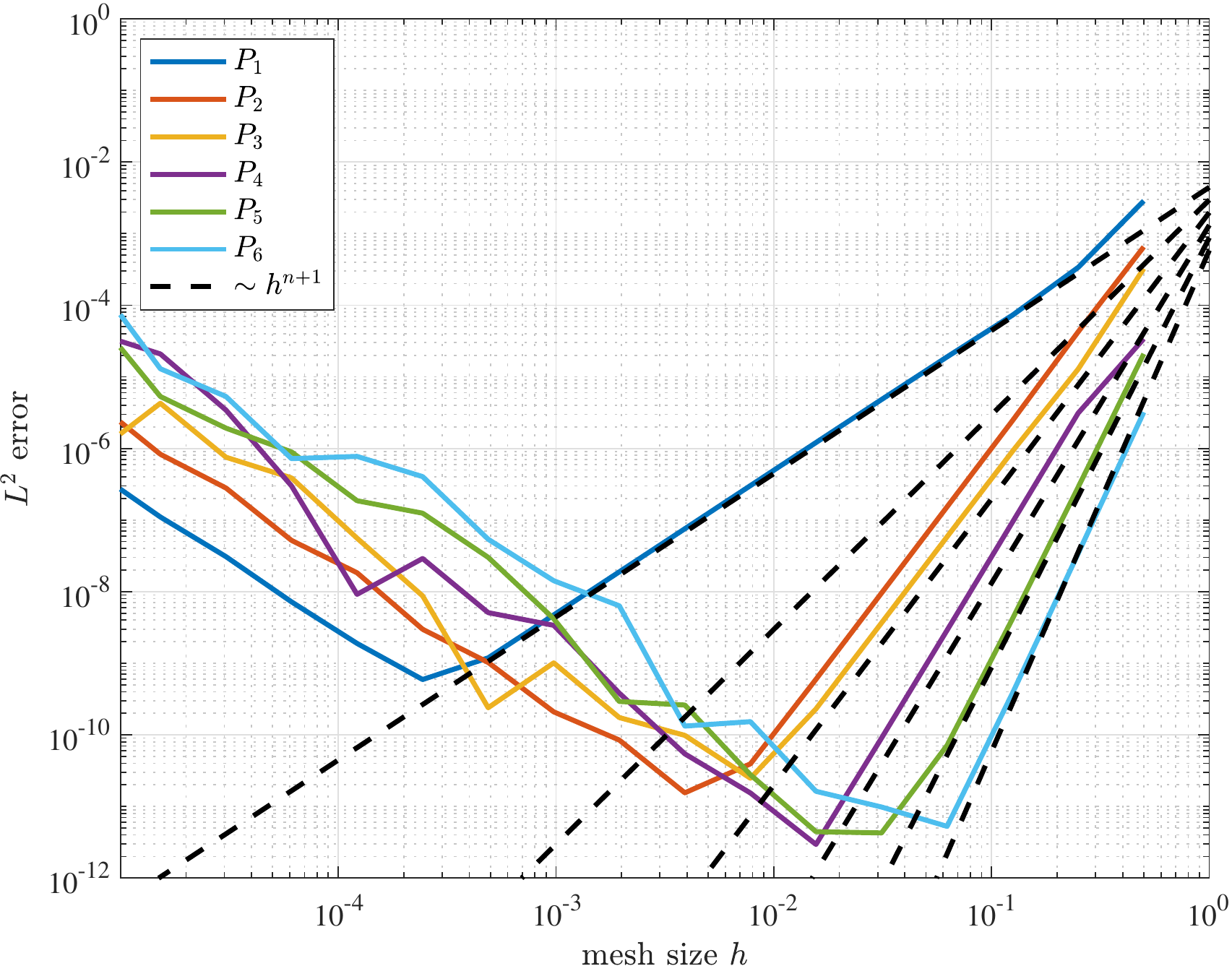}
  \caption{$L^2$ experimental convergence order of solution to the problem $Lu=f$ with $f(x)=x$ for different FE Lagrange elements P$_n$ (left) with degenerate $\mu=x^2$ remains first order and (right) with nondegenerate $\mu=1+x^2$ shows optimal convergence order. The error goes down to $\sim 10^{-11}$, similar to finite difference methods, it goes up due to ``loss of significance''.}
  \label{fig:suboptimal_convergence}
\end{figure}
}

\bibliographystyle{ieeetr}
\bibliography{thinfilm} 
\end{document}